\documentclass[a4paper,12pt]{article}
\usepackage{amssymb}
\usepackage{amsthm}
\usepackage{amsmath}
\usepackage{amsfonts}

\begin{document}

\def\sect{\section}

\newtheorem{thm}{Theorem}[section]
\newtheorem{cor}[thm]{Corollary}
\newtheorem{lem}[thm]{Lemma}
\newtheorem{prop}[thm]{Proposition}
\newtheorem{propconstr}[thm]{Proposition-Construction}

\theoremstyle{definition}
\newtheorem{para}[thm]{}
\newtheorem{ax}[thm]{Axiom}
\newtheorem{conj}[thm]{Conjecture}
\newtheorem{defn}[thm]{Definition}
\newtheorem{notation}[thm]{Notation}
\newtheorem{rem}[thm]{Remark}
\newtheorem{remark}[thm]{Remark}
\newtheorem{question}[thm]{Question}
\newtheorem{example}[thm]{Example}
\newtheorem{problem}[thm]{Problem}
\newtheorem{excercise}[thm]{Exercise}
\newtheorem{ex}[thm]{Exercise}

\def\Bbb{\mathbb}
\def\cal{\mathcal}
\def\mL{{\mathcal L}}
\def\mC{{\mathcal C}}

\overfullrule=0pt

%%Macros

%\def\Bbb{\bf}

\def\si{\sigma}
\def\prf{\smallskip\noindent{\it	Proof}.	}
\def\call{{\cal	L}}
\def\nat{{\Bbb	N}}
\def\la{\langle}
\def\ra{\rangle}
\def\inv{^{-1}}
\def\trdeg{{tr.deg}}
\def\dim{{\rm	dim}}
\def\th{{\rm	Th}}
\def\rest{{\lower	.25	em	\hbox{$\vert$}}}
\def\ch{{\rm	char}}
\def\zee{{\Bbb	Z}}
\def\conc{^\frown}
\def\cals{{\cal	S}}
\def\mult{{\rm	Mult}}
\def\calv{{\cal	V}}
\def\aut{{\rm	Aut}}
\def\ffi{{\Bbb	F}}
\def\ffiti{\tilde{\Bbb		F}}
\def\calg{{\cal G}}
\def\gal{{\cal G}al}
\def\calc{{\cal C}}
\def\rat{{\Bbb Q}}
\def\ione{{\it (I1)}}
\def\itwo{{\it (I2)}}
\def\itwow{{\it (I2w)}}
\def\ithreew{{\it (I3w)}}
\def\ithree{{\it (I3)}}
\def\ifour{{\it (I4)}}
\def\ifourw{{\it (I4w)}}
\def\bigudot{\hbox{$\bigcup\kern-.75em\cdot\;\,$}}
\def\caln{{\cal N}}
\def\cald{{\cal D}}
\def\calm{{\cal M}}
\def\calu{{\cal U}}
\def\Div{\hbox{Div}}
\def\pic{\hbox{Pic}}

\def\G{{\cal G}}
\def\sime{{\hbox{$\kern-.19em\sim$}}}
\def\calp{{\cal P}}
\def\asi{almost-$\cals$-internal}
\def\api{almost-$\{p\}$-internal}
\def\acb{\overline{\rm Cb}}
\def\dnfo{\,\raise.2em\hbox{$\,\mathrel|\kern-.9em\lower.35em\hbox{$\smile$}
$}}
\def\dnf#1{\lower1em\hbox{$\buildrel\dnfo\over{\scriptstyle #1}$}}
\def\dfo{\;\raise.2em\hbox{$\mathrel|\kern-.9em\lower.35em\hbox{$\smile$}
\kern-.7em\hbox{\char'57}$}\;}
\def\df#1{\lower1em\hbox{$\buildrel\dfo\over{\scriptstyle #1}$}}

\def\cl{{\rm cl}}
\def\ho{hereditarily orthogonal\ }
\def\acl{{\rm acl}}
\def\calt{{\cal T}}
\def\Fix{{\rm Fix }}
\def\dcl{{\rm dcl}}
\def\fix{{\rm Fix}}
\def\frob{{\rm Frob}}
\def\Cb{{\rm Cb}}
\def\Frob{{\rm Frob}}
\def\vlabel{\label}
\def\Th{{\rm Th}}
\def\pia#1{$\{#1\}$-analysable}
\def\Stab{{\rm Stab}}
\def\ld{{\rm ld}}
\def\ild{{\rm ild}}

\title{A note on canonical bases and one-based types in
supersimple
theories}
\author{Zo\'e Chatzidakis\thanks{A revised version of this
paper was written at the Newton Institute during the special semester 
``Model Theory and Applications to Algebra and Analysis'', and the
author gratefully aknowledges their support. The author was also
partially supported by  MRTN-CT-2004-512234 and by ANR-06-BLAN-0183.}
(IMJ, CNRS - Universit\'e
    Paris 7)} 
%% \centerline{}
%% \centerline{Zo\'e Chatzidakis\footnote{$^{*}$}{A revised version of this
%% paper was written at the Newton Institute during the special semester 
%% ``Model Theory and Applications to Algebra and Analysis'', and the
%% author gratefully aknowledges their support. The author was also
%% partially supported by  MRTN-CT-2004-512234 and by ANR-06-BLAN-0183}  --
%% CNRS (Universit\'e 
%% Paris 7)} 

\maketitle
%\bigskip\noindent
\begin{abstract}This paper studies the CBP, a model-theoretic property
  first discovered 
  by Pillay and Ziegler. We first show a general decomposition result
  of types of canonical bases, which one can think of as a sort of primary
  decomposition. 
This decomposition is then used to show that existentially
  closed difference fields of any characteristic have the CBP. 
 We also derive consequences of the CBP, and use these results for
 applications to differential and difference varieties, and algebraic
 dynamics.
\end{abstract}
\section*{\bf Introduction}

In \cite{P}, Anand Pillay gives a model-theoretic translation of a
property enjoyed by compact complex manifolds (and proved 
by Fr\'ed\'eric Campana and by Akira Fujiki). With Martin Ziegler,  he then shows in
\cite{[PZ]} that various algebraic structures enjoy this property
(differentially closed fields of characteristic $0$; existentially
closed difference fields of characteristic $0$). As 
with compact complex manifolds, their proof has as immediate consequence 
the dichotomy for types of rank $1$ in these algebraic structures. 
This property will later be called the  {\em Canonical Base Property}  (CBP for
short) by Rahim Moosa and 
Pillay \cite{[MP]}.  We will state the precise
definition of the CBP later (see \ref{cb0}), as it 
requires several model-theoretic definitions, but here is a rough idea. Let us assume that we have good notions of
independence, genericity and dimension, and let $S\subset X\times Y$ be
definable. Viewing $S$ as a family of definable subsets $S_x$ of $Y$,
assume that for $x\neq x'$ in $X$, $S_x$ and $S_{x'}$ do not have the
same generics, and have finite dimension. Fix some $a\in X$, a generic
$b$ of $S_a$. The CBP then gives strong restrictions on the set
$S^b=\{x\in X\mid b\in S_x\}$: for instance in the complex manifold case, it is
Moishezon, and in the differential case it is isoconstant.

The aim of this paper is three-fold: give
reductions to prove the CBP; derive consequences of the CBP; show that
existentially closed difference fields of positive characteristic have
the CBP. We then give some applications of these results to 
differential  and difference varieties.   

We postpone a
detailed description of the model-theoretic results of this paper to the
middle of section 1 (\ref{cb1}) and to the beginning of section 2, but we will now describe two of the
algebraic applications. First, an algebraic consequence of Theorems
\ref{thm1} and \ref{prop2}. 
We work in some large existentially closed 
difference field 
$(\calu,\si)$, of characteristic $p$; if $p>0$, $\frob$ denotes the map
$x\mapsto x^p$ and if $p=0$, the identity map. 

\medskip\noindent
{\bf Theorem \ref{thm1}$'$}. Let $A,B$ be  difference subfields 
of $\calu$ intersecting in $C$, with algebraic closures intersecting
in $C^{alg}$, and with $\trdeg(A/C)<\infty$. Let $D\subset B$ be generated
over $C$ by all tuples $d$ such that there exist an algebraically
closed difference field  $F$ containing $C$ and free from $B$ over $C$,
and integers $n>0$ and $m$ such that $d\in F(e)$ for some tuple $e$
of elements satisfying $\si^n\frob^m(x)=x$. Then $A$ and $B$ are free
over $D$. 

\medskip
The  purely model theoretic result
\ref{descent} yields descent results for differential
and difference varieties (\ref{descent1} and \ref{descent2}). We state
here a consequence in terms of algebraic dynamics:

\medskip\noindent
{\bf Theorem \ref{descent3}}. Let $K_1,K_2$ be fields
intersecting in $k$ and with algebraic closures intersecting in $k^{alg}$; for $i=1,2$, let $V_i$ be an absolutely irreducible
variety and 
$\phi_i:V_i\to V_i$ a dominant rational map defined over $K_i$. Assume  that $K_2$ is a regular extension of $k$ and 
that there are an integer $r\geq 1$ and a dominant rational map $f:V_1\to V_2$ 
such that $f\circ \phi_1=\phi_2^{(r)} \circ f$ (where $\phi_2^{(r)}$
denotes the map obtained by iterating $r$ times $\phi_2$). Then there is
a variety $V_0$ and a 
dominant rational map $\phi_0:V_0\to V_0$, all defined over $k$, a dominant
map $g:V_2\to V_0$ such that $g\circ \phi_2=\phi_0\circ   g$, and
$\deg(\phi_0)=\deg(\phi_2)$.  

\medskip
The particular way this result is stated is motivated by a question of
Lucien 
Szpiro and Thomas Tucker concerning descent for algebraic
dynamics, arising out of Northcott's theorem for dynamics over function
fields. Assume that $K_2$ is a function field over $k$,
and that some {\em limited}\footnote{see \cite{AD1} for a definition} subset $S$
of $V_2(K_2)$ satisfies that  $\bigcap_{j=0}^n \phi^{(j)}(S)$ is Zariski
dense in $V_2$ for every $n>0$. One
can then find $(V_1,\phi_1)$, $r$ and $f$ as above, so that our result
applies to 
give a quotient $(V_0,\phi_0)$ of $(V_2,\phi_2)$ defined over the
smaller field $k$ and with $\deg(\phi_0)=\deg(\phi_2)$. Under certain
hypotheses, one can even have this $g$ be birational, see \cite{AD1},
\cite{AD2}.

\medskip
This note originally contained a proof that a
type analysable in terms
of one-based types is one-based. However, Frank Wagner \cite{[W2]}  found a
much nicer
proof, working in a more general context, so that this part of the
note disappeared.
A result of independent interest, Proposition~\ref{cor1}, obtained as a
by-product of the study of one-based types and appearing in the
appendix, tells us that if $p$ is a type of SU-rank $\omega^\alpha$ for
some ordinal $\alpha$ and with algebraically
closed base of finite SU-rank, then there is a smallest algebraically
closed set over 
which there is a type of SU-rank $\omega^\alpha$ non-orthogonal to
$p$. The condition of finite 
rank of the base is necessary.

\bigskip\noindent
The paper is organised as follows. Section 1 contains all definitions
and preliminary results on supersimple theories, as well as the proof of
the decomposition result 
Theorem \ref{thm37}. Section~2 contains  various results which are
consequences of 
the CBP. Section 3 shows that if $K$ is an existentially closed
difference field of any characteristic, then $\Th(K)$ has the
CBP. Section 4 contains some applications of the CBP to differential  and
difference varieties. Section 5 is the appendix. 

\smallskip
Some words on the chronology of the paper and results on the CBP. It all
started with 
the result of Pillay and Ziegler \cite{[PZ]}, a result
inspired by a result of Campana on compact complex spaces (see \cite{P}), and which
prompted me to look
at the general case. The first version of this paper, which
contained only Theorem~\ref{thm1}, an old version of Theorem~\ref{prop2} and
Proposition~\ref{cor1}, as well as the proof that a type analysable in
one-based types was one-based, was written in 2002. Almost instantly the
result on analysable one-based types was generalised by Wagner. The
paper was submitted, but not accepted during several years. In the
meantime, Moosa and Pillay, having read and believed the preprint,
further investigated the CBP in \cite{[MP]}. Reading their preprint
alerted me to 
the fact that the CBP might imply other stronger properties, as
suggested by the fact that compact complex analytic spaces had the
UCBP. Thus the material in section 2 starting  from \ref{lem2-1} on,
came later (end of 2008, and 2011). Independently,  Prerna
Juhlin (\cite{J}) has 
obtained several results 
on theories with the CBP in her doctoral thesis (2010). Moosa studies in
\cite{Mo} variants of internality in presence of the CBP. Daniel Palacin and
Wagner  continue and generalise the study of the CBP in
\cite{PW}. Ehud Hrushovski (\cite{H}) gives an example of an
$\aleph_1$-categorical theory which does not have the CBP. This example
now appears in a paper by Hrushovski, Daniel Palac\'\i n and Pillay
\cite{HPP}.

\sect{Results on supersimple theories}
\para\vlabel{sett}{\bf Setting}. We work in a model $M$ (sufficiently
saturated)
of
a  complete theory $T$, which is
supersimple and eliminates imaginaries.  The
results given below 
generalise easily to a simple theory eliminating hyperimaginaries,
provided that some of the sets considered are ranked by the SU-rank.

Given (maybe infinite) tuples $a,b\in M$, we  denote by  $\acb(a/b)$
the smallest algebraically closed
subset of $M$ over which $tp(a/b)$ does not fork.  Since our
theory is supersimple, it coincides with the algebraic
closure (in $M^{eq}$) of the usual canonical basis $\Cb(a/b)$ of
$tp(a/b)$, and
is contained in $\acl(b)$.
For classical results on canonical bases and supersimple theories,
    see e.g. sections 3.3 and 5.1 -- 5.3 of \cite{[W1]}.
We will use repeatedly the following
consequences of
our hypotheses on $T$:

\begin{enumerate}
\item{Let $ B\subset M$, $a\in M$, and $(a_n)_{n\in\nat}$ a
sequence of $B$-independent realizations of $tp(a/B)$. Then for some
$m$, $\acb(a/B)$ is contained in $\acl(a_1\ldots a_m)$; for any $n$, 
$\acl(a_1\ldots a_n)\cap B\subseteq \acb(a/B)$.}

\item{Let $B\subset M$, $a\in M$, and $(a_n)_{n\in\nat}$ a
sequence of $B$-independent realizations of $tp(a/B)$. Let $m$ be minimal such
that  
$C=\acb(a/B)\subseteq \acl(a_1\ldots a_m)$. Then
$SU(a_1/a_2\ldots a_m)>SU(a/B)$: otherwise
$a_1\dnfo_{a_2\ldots a_m}C$ would imply $C\subseteq \acl(a_2\ldots
a_m)$, and contradict the minimality of $m$.}
   \item{If $A$ and $B$ are algebraically
   closed subsets of $M$ intersecting in $C$, and $D$ is independent
   from
   $AB$ over $C$, then $\acl(DA)\cap \acl(DB)=\acl(DC)$ (if $e\in
   \acl(DA)\cap \acl(DB)$, then $\acb(De/AB)\subseteq A\cap B=C$).}
\end{enumerate}

\para{\bf Internality and analysability}.
In what follows, we will assume that $\cals$ is a set of types with
algebraically closed base and which
       is closed under  $\aut(M/\acl(\emptyset))$-conjugation. Then non-orthogonality generates  an
equivalence relation on the regular types in $\cals$. For more details,
see section 3.4 of \cite{[W1]}. 

Recall that if  $a\in M$ and $A\subseteq M$, then $tp(a/A)$
is {\em $\cals$-internal} [resp., {\em almost-$\cals$-internal}] if
there is
some set $B=\acl(B)$
containing
$A$ and independent from $a$ over $A$,  and a  tuple
$b_1,\ldots,b_n$ such that $a\in \dcl(Bb_1\ldots b_n)$ [resp.
$a\in
\acl(Bb_1\ldots b_n)$], and each $b_i$ realises a type which is in
$\cals$ and has base contained in $B$. % Note that if $tp(a/A)$ is
% \asi, then $a$ is equi-algebraic
% over $A$ with a tuple whose type over $A$ is $\cals$-internal:
% consider $\Cb(Bb_1\ldots b_n/\acl(Aa))$.

$tp(a/A)$ is {\em $\cals$-analysable} if there are $a_1,\ldots,a_n$
such
that $\acl(Aa_1\ldots a_n)=\acl(Aa)$ and each
$tp(a_i/Aa_1\ldots a_{i-1})$ is $\cals$-internal (or equivalently,
each
$tp(a_i/Aa_1\ldots a_{i-1})$ is \asi).

\para\vlabel{obs}{\bf Observations}. Let  $A=\acl(A)\subset M$.
\begin{enumerate}
\item{If $tp(a_i/A)$ is almost-$\cals$-internal for
$i=1,\ldots,n$, then so is $tp(a_1\ldots a_n/A)$.}
\item{If $tp(a/A)$ is \asi, and $b\in \acl(Aa)$, then $tp(b/A)$
is \asi.}
\item{If $\cals'$ is a set of types which are \asi,
and if $p$ is almost-$\cals'$-internal, then $p$ is \asi. }
\item{Similarly for
$\cals$-analysability.}
\item{Let $\cals_1$, $\cals_2\subset \cals$ be sets of types of SU-rank $1$
which are closed under $\aut(M)$-conjugation. 
If all types in $\cals_1$ are orthogonal to all types in $\cals_2$
(denoted by $\cals_1\perp\cals_2$) and  
  $q_i$ is $\cals_i$-analysable for $i=1,2$, then all
extensions of $q_1$ 
are orthogonal to all extensions of $q_2$.}

\end{enumerate}

\para\vlabel{mod1}{\bf One-basedness}. Let $S\subseteq M^k$ be
$A$-invariant.
Then $S$ is {\em one-based (over $A$)} if whenever $b$ is a tuple of
elements of $S$, and $B\supseteq A$ then $b$ is independent from $B$
over $\acl(Ab)\cap \acl(B)$. A type $p$ (over $A$) is  {\em one-based}
if the set of its
realisations is one-based over $A$. 

\smallskip\noindent
{\bf Properties} (see \cite{[W2]}). 
\begin{enumerate}
\item{Let $p$ be a type, and $q$ a non-forking extension of
$p$. Then $p$ is one-based if and only if $q$ is
one-based. One-basedness is preserved under $\aut(M)$-conjugation.} 
\item{A union of one-based sets is one-based.}
\item{A type analysable by one-based types is one-based.}
\end{enumerate}

\para\vlabel{nonorth}{\bf Non-orthogonality and internality}.
 Let $A=\acl(A)\subset M$, and $a$ a tuple
in $M$, with 
$SU(a/A)=\beta+\omega^\alpha$ for some $\alpha$, $\beta$. Then there is
$B=\acl(B)$ containing $A$ and independent from $a$ over $A$, and $b$
such that $SU(b/B)=\omega^\alpha$ and $b\dfo_B a$ (see \cite{[W1]}, 5.1.12). Then
$C=\acb(Bb/Aa)$ is contained in the algebraic closure of independent
realisations of $tp(Bb/\acl(Aa))$, and therefore its type over $A$ is
almost internal to the set of conjugates of $tp(b/B)$ over $A$. 

If $tp(b/B)$ is one-based, then $tp(C/A)$ is one-based (by
\ref{mod1}(3)), so that 
$\acl(Bb)$ and $C$ are independent over their intersection $D$, and
therefore $C=D$. Since
$SU(b/B)=\omega^\alpha$, a standard computation gives $SU(C/A)=\omega^\alpha$.

\para\vlabel{semimin}{\bf Semi-minimal analysis}.
In particular every finite
SU-rank type has a {\em semi-minimal analysis}, i.e.: given
$A=\acl(A)$ and
$a$ of finite SU-rank over $A$, there are tuples $a_1,\ldots,a_n$ such
that
$\acl(Aa_1 \ldots a_n)=\acl(Aa)$, and for every
$i$, either $tp(a_i/Aa_1\ldots a_{i-1})$ is one-based of SU-rank $1$,
or
it is internal to the set of conjugates of some non-one-based
type of
SU-rank $1$.

\para\vlabel{cb0}{\bf Definition of the CBP}. Let $T$ be a simple
theory, which eliminates imaginaries and hyperimaginaries.  The theory
$T$ {\em has the  
 CBP} if 
whenever $A$ and $B$ are algebraically closed sets of finite SU-rank
over their intersection, and $A=\acb(B/A)$, then $tp(A/B)$ is \asi{},  where
$\cals$  is the set of types of SU-rank $1$ with algebraically closed
base. [Actually, as we will see in Theorem \ref{thm37}, it suffices to
take  for $\cals$ the set of non-one-based types of SU-rank $1$ with
algebraically closed base.]\\[0.1in]
One can also restrict this definition to smaller families of types: let 
$\calp$ be a family of types of finite SU-rank and with algebraically
closed base. We 
say that $\calp$ {\em has the CBP} if whenever $D$ is algebraically closed,  $b$ is a tuple of realisations
of types in $\calp$ with base contained in $D$, 
and  $A=\acb(Db/AD)$, then $tp(A/\acl(Db))$ is \asi, for the family
$\cals\subset \calp$ of types in $\calp$ of SU-rank $1$ [and which are
not one-based]. Thus Pillay and Ziegler
show in \cite{[PZ]} that the family of very thin types in separably closed fields of
finite degree of imperfection has the CBP.  See  the concluding remarks
at the end of 
section 2 for a discussion.

\para{\bf Definition}\vlabel{ho}. Let $p$ and $q$ be types. We say that
$p$ is   {\em hereditarily 
orthogonal} to $q$ if every extension of $p$ is orthogonal to $q$.

\para\vlabel{lem1}{\bf Lemma}. Let $E,B\subset M$
be algebraically closed sets,
$b\in M$ a tuple.  Assume that $tp(b/B)$ is
\asi, $E=\acb(Bb/E)$, and $\cals$ is closed under $\aut(M/E)$-conjugation.
If $A=\acb(B/E)$, then $tp(E/A)$ is \asi.

\prf Let $(B_1b_1),\ldots,(B_nb_n)$ be realisations of
$tp(Bb/E)$ which are independent  over $E$ and such that  $E\subseteq
\acl(B_1b_1\ldots B_nb_n)$. Since $B\dnfo_AE$, we get 
$B_1\ldots B_n\dnfo_AE$; observation \ref{obs}(3) then   gives the result.

\para\vlabel{lem2}{\bf Lemma}. Let $E, F\subset M$ be algebraically
closed sets, with $\acb(E/F)=F$. If $E_0=\acb(F/E)$, then
$F=\acb(E_0/F)$.

\prf Let $F_0=\acb(E_0/F)$. Then $E_0\dnfo_{F_0}F$ and $E\dnfo_{E_0}F$, which imply
$E\dnfo_{E_0F_0}F$ (since $F_0\subseteq F$) and $E\dnfo_{F_0}F$ by
transitivity. Hence $F_0=F$.

\para\vlabel{lem36}{\bf Lemma}. Let $\cals_1$ and $\cals_2$ be sets
of types of SU-rank $1$  closed under 
$\aut(M)$-conjugation, with $\cals_1\perp \cals_2$.
Assume
that $tp(E_i)$ is $\cals_i$-analysable for $i=1,2$, and that
$D=\acl(D)\subseteq \acl(E_1E_2)$. Let $D_i=D\cap \acl(E_i)$ for $i=1,2$.
Then
$D=\acl(D_1D_2)$.

\prf Without loss of generality, each $E_i$ is algebraically closed.
Since $\acb(E_1/D)$ realises an $\cals_1$-analysable type, it equals
$D_1$ and hence $D\dnfo_{D_1}E_1$. As $D\subseteq \acl(E_1E_2)$, this
implies that $tp(D/D_1)$ is $\cals_2$-analysable. Hence so is
$tp(D/D_1D_2)$.
Similarly, $D\dnfo_{D_2}E_2$ and $tp(D/D_1D_2)$ is $\cals_1$-analysable.
Our hypothesis on the orthogonality of the members of $\cals_1$ and
those of $\cals_2$ then implies $D\subseteq
\acl(D_1D_2)$: a type which is $\cals_1$-analysable and
$\cals_2$-analysable must be algebraic. 

\para\vlabel{lem25}{\bf Lemma}. Let $\cals$ be a set of types of
SU-rank $1$, which is closed under $\aut(M)$-conjugation. Let $B\subset
F$ and $A$ be algebraically closed 
sets such that $tp(A)$  and $tp(B)$ are
\asi\ [resp. $\cals$-analyzable], and $B$ is maximal 
contained in $F$ with
this property. 
\begin{enumerate}
\item{Then $\acl(AB)$ is the maximal subset of $\acl(AF)$ whose 
 type is \asi{} [resp. $\cals$-analyzable].}
\item{Let $G$ be independent from $F$. Then $\acl(GB)$ is the
maximal subset of $\acl(GF)$ whose type over $G$ is \asi\ [resp. $\cals$-analyzable].}
\end{enumerate}

\prf (1) Let $d\in \acl(AF)$ be such that $tp(d)$ is \asi. Then so is the
type of $\acb(Ad/F)$; hence $\acb(Ad/F)\subseteq B$ and $d\in
\acl(AB)$. Same proof for  $\cals$-analyzable. 

(2) Let $e\in\acl(GF)$ realise an \asi\ type over $G$. As $G\dnfo F$,
    $\acb(Ge/F)$ realises an \asi\ type, hence is contained in
     $B$. Hence $Ge\dnfo_BF$, which implies $e\in\acl(GB)$. Same proof
for  $\cals$-analyzable.

\medskip The following result is well-known, but for lack of a
reference, we will give the proof. 

\para\vlabel{lem41}{\bf Lemma}. Let $p$ and $q$ be types over sets $A$
and $B$ respectively, and assume that  $p\not\perp q$. Then for some
integer $\ell$ there are
realisations $a_0,\ldots,a_\ell$ of $p$, $b_0,\ldots,b_\ell$ of $q$, such that
the  tuples $a_i$ are independent over $A$, the tuples $b_j$ are independent
over $B$, $$a_0\ldots a_\ell\dnfo_AB, \ b_0\ldots b_\ell\dnfo_BA \hbox{ and
 } a_0\ldots a_\ell\dfo_{AB}b_0\ldots b_\ell.$$ 

\prf 
By assumption there are some $C$ containing $A$ and $B$, and
realisations $a$ of $p$, $b$ of $q$ such that $a\dnfo_AC$, $b\dnfo_BC$
and $a\dfo_Cb$. Let $D=\acb(a,b/C)$. Then for some $\ell$ there are independent
realisations $(a_i,b_i)$, $i=1,\ldots,\ell$, of $tp(a,b/C)$ such that $D\subset
\acl(ABa_1\ldots a_\ell b_1\ldots b_\ell)$ (by \ref{sett}(1)); we may
choose these realisations to be  independent from $(a,b):=(a_0,b_0)$ over $C$. 
Then $$a_0\dfo_{ABa_1\ldots
a_\ell b_1\ldots b_\ell}b_0.$$ As
$a\dnfo_AC$, $b\dnfo_BC$ and $C$ contains $AB$, the tuples $a_i$ and $b_j$
also satisfy 
the required first four conditions. Transitivity of independence then
implies
$$a_0\ldots a_\ell\dfo_{AB}b_0\ldots b_\ell.$$ 

\para\vlabel{lem35}{\bf Notation and some remarks}. Let $p$ be a type
of SU-rank $1$ over some
algebraically closed set, and let $C=\acl(C)$.
We denote by $\cals(p,C)$ the smallest set of types of SU-rank $1$
with algebraically closed base,
which contains  $p$
and is closed under
$\aut(M/C)$-conjugation.
We write $\cals(p)$ for
$\cals(p,\emptyset)$.

Let $p$ and $q$ be types of SU-rank $1$, with algebraically closed base
$A$ and $B$ respectively. Certainly if $q$ is \api, then $p\not\perp
q$. If $A=B$, then the converse holds: $p\not\perp q$ iff $q$ is
\api\ (iff $p$ is almost-$\{q\}$-internal). If $A\neq B$, then
$p\not\perp q$ implies that $q$ is almost-$\cals(p,B)$-internal, since
any two realisations of $q$ are in the 
same 
$\aut(M/B)$-orbit; but in general, $q$ will not be \api.   

In particular, if
$q\not\perp p$, then $q$ is almost-$\cals(p)$-internal. Hence, 
\begin{center}either
$\cals(p)\perp \cals(q)$,
or every member of $\cals(p)$ is $\cals(q)$-internal\end{center}(and every member of
$\cals(q)$ is $\cals(p)$-internal).

\bigskip
In the rest of the first two sections of the paper, the letters $\cals$,
$\cals'$, $\cals_1$, 
etc. will always denote sets of SU-rank $1$ types with algebraically
closed base.   

\para\vlabel{cb1}
We now start towards the proof of Theorem \ref{thm37}. It will reduce the problem of showing the CBP to
showing it for $\{p\}$-analysable types when $p$ is a type of SU-rank
$1$ with algebraically closed base. This reduction is essential in the proof that existentially closed difference fields 
of positive characteristic have the CBP. We conclude the section with
small partial results.

\para\vlabel{lem37}{\bf Proposition}.  Let $F$ and $E$ be algebraically
closed sets such that $F\cap
E=C$, $SU(E/C)$ and $SU(F/C)$ are finite, $\acb(E/F)=F$ and
$\acb(F/E)=E$.  There are non one-based types $p_1,\ldots,p_m$ of
SU-rank $1$, algebraically closed sets 
$E_1,\ldots,E_m, F_1,\ldots, F_m$ such that, letting
$\cals_i=\cals(p_i,C)$ for $i=1,\ldots,m$,
\begin{itemize}
\item[(i)]  $tp(E_i/C)$ and $tp(F_i/C)$ are $\cals_i$-analysable, 
$\acb(E_i/F_i)=F_i$ and $\acb(F_i/E_i)=E_i$, and
\item[(ii)] $\acl(E_1\ldots E_m)=E$, $\acl(F_1\ldots F_m)=F$.
\item[(iii)] The sets $E_i$ are independent over $C$, as well as the
sets  $F_i$. 
\end{itemize}

\prf Assume the
result false, and take a counterexample with $SU(EF/C)$ minimal among
all possible $(E,F,C)$, and among those, with $SU(F/C)+SU(E/C)$ minimal.

Let $p_1,\ldots,p_m$ be types of SU-rank $1$ with algebraically closed
base, which are pairwise orthogonal,  such that each
$p_i$ is non-orthogonal to $tp(E/C)$ or to $tp(F/C)$, and such that
any SU-rank $1$ type which is non-orthogonal to one of $tp(E/C)$,
$tp(F/C)$, is non-orthogonal to one of the types $p_i$ (see section 5.2
in \cite{[W1]}). We let $\cals_i=\cals(p_i,C)$, and   $E_i$ and $F_i$ the
maximal subsets of $E$ and $F$ respectively such that  $tp(E_i/C)$ and
$tp(F_i/C)$ 
are $\cals_i$-analysable. We need to show that
no $p_i$ is one-based, the second part of item (i), and item (ii) (item
(iii) is immediate since the types $p_i$ are pairwise orthogonal, by
  \ref{obs}(5)). 

\smallskip
Adding to the
language constants symbols for the elements of $C$,
{\em we will assume that $C=\emptyset$}.

\smallskip
We  say that a set $D$ {\em satisfies $(*)$ over a set $H$} if
$D=\acl(HD_1\ldots D_m)$, where $tp(D_i/H)$ is $\cals_i$-analysable for
each $i$. 
If $H=\acl(\emptyset)$, we will
simply say that $D$ satisfies $(*)$. Note that by \ref{lem36}, any
subset of $\acl(HD)$ whose type over $H$ is $\cals_i$-analysable will be
contained in $\acl(HD_i)$. 

By Lemma \ref{lem36},
an algebraically closed subset of a set satisfying $(*)$ over $H$ also
satisfies $(*)$ over $H$, and  the algebraic closure of a union of sets
satisfying $(*)$ over $H$ satisfies $(*)$ over $H$.  Hence, if $D$
satisfies $(*)$ over $H$, $J\supseteq H$, and
$\acb(D/J)=J$, then $J$ satisfies $(*)$ over $H$, as $J$ is contained in the
algebraic closure of finitely many realisations of $tp(D/H)$.

Assume that   $E$  satisfies
$(*)$; then $F=\acb(E/F)$
 satisfies $(*)$, and therefore can be written as $\acl(F'_1,\ldots,F'_m)$,
 where each $tp(F'_i)$ is $\cals_i$-analysable. Each $F'_i$ is contained
 in $F_i$, and therefore $\acl(F_1 \ldots F_m)=F$ and $F_i=F'_i$. We
 know that $F$ is 
 contained in the algebraic closure of $F$-independent realisations of
 $tp(E/F)$. Lemma \ref{lem36} then
 gives us that necessarily $F_i$ is contained in the algebraic closure
 of $F$-independent realisations of $tp(E_i/F)$. Then \ref{sett}(1)
 implies that  
 $\acb(E_i/F)\supseteq F_i$; the reverse inclusion holds since
 $tp(\acb(E_i/F))$ is $\cals_i$-analysable. By symmetry,
$E_i=\acb(F_i/E)$. 
  Furthermore, no $p_i$ is
one-based: otherwise, by \ref{mod1}(3) $tp(E_i)$ would be
one-based,
whence $E\cap F=\acl(\emptyset)$ would yield $E_i\dnfo F$, and therefore
$E_i=F_i=\acl(\emptyset)$. 
This shows that if $E$ satisfies $(*)$, then the conclusion
of the Lemma holds. By symmetry  neither $E$ nor $F$ satisfies $(*)$.

Using the semi-minimal analysis of $tp(F)$, there is
$B=\acl(B)\subset F$, $B\neq F$,  and a type $p$ of SU-rank $1$, such
that
$tp(F/B)$ is
almost-$\cals(p)$-internal. Note that $B\neq \acl(\emptyset)$:
otherwise $tp(F/C)$ would be $\cals(p)$-internal, contradicting our
assumption that $F$ does not satisfy $(*)$. 
Let $A=\acb(B/E)$. Then $tp(E/A)$
is almost-$\cals(p)$-internal
by Lemma \ref{lem1}, so that $A\neq \acl(\emptyset)$.

\smallskip\noindent
{\bf Step 1}. $A$ satisfies  $(*)$.

Let  $B_0=\acb(A/B)$. Then
$B_0\neq \acl(\emptyset)$ (because
otherwise $B$ and $A$ would be independent), and $\acb(B_0/A)=A$ by
Lemma \ref{lem2}.
As $SU(B_0)<SU(F)$ and $\acl(AB_0)\subseteq \acl(EF)$, by
      induction hypothesis  $A$ satisfies $(*)$.

\smallskip
Thus $A\neq E$. Since $tp(E/A)$ is almost-$\cals(p)$-internal, if
$B_1=\acb(A/F)$, then $tp(F/B_1)$ is
almost-$\cals(p)$-internal by Lemma \ref{lem1},   $B_1$
satisfies $(*)$, and $\acl(\emptyset)\neq B_1\neq F$. 
Let
$E_p$ be the largest
subset of $E$ realising an $\cals(p)$-analysable type. If $p$ is orthogonal
to every $p_i$, then $E_p=\acl(\emptyset)$, by definition of the set
$\{p_1,\ldots, p_m\}$. Otherwise, by \ref{lem35}, $\cals(p)=\cals_i$ for
some $i$, and therefore $E_p=E_i$; we will
first show in the next two steps that this case is impossible.

\smallskip\noindent
{\bf Step 2}. $\acl(FE_p)\cap E=E_p$.

Let $D=\acl(FE_p)\cap E$. If $D\neq E_p$, then, using the semi-minimal
analysis of $tp(D/E_p)$, there is $d\in D\setminus E_p$ with $tp(d/E_p)$
almost-$\cals(q)$-internal for some  type  $q$  of SU-rank
$1$. By maximality of 
$E_p$, we have $\cals(q)\perp \cals(p)$. Since $tp(F/B_1)$ is
almost-$\cals(p)$-internal, we obtain $d\in \acl(B_1E_p)$.
Because $B_1$ and $E_p$ satisfy $(*)$, and $tp(d/E_p)\perp p$,  using \ref{lem36} we may write $\acl(E_pd)$ as
$\acl(E_pD_1)$, where  $tp(D_1)$ is
$\cals(q)$-analysable and 
hereditarily
orthogonal to $p$. Thus $D_1\dnfo E_p$; because $D_1\subseteq
\acl(B_1E_p)$, we obtain $D_1\subseteq B_1$; as $D_1\subseteq E$, this
implies $D_1=
\acl(\emptyset)$,  a contradiction.

\smallskip\noindent
{\bf Step 3}. $E_p=\acl(\emptyset)$.

Let $D=\acb(E/FE_p)$. Then $\acb(E/D)=D$. By Step 2, $D\cap E=E_p$.
Moreover $\acb(D/E)=E$:
let $E_0=\acb(D/E)$; from $E\dnfo_DF$ we deduce $E\dnfo_{DE_0}F$; since
$E\dnfo_{E_0}D$,  transitivity gives $E\dnfo_{E_0}F$, and therefore $E=E_0$.
Thus, if $E_p\neq
\acl(\emptyset)$, then $SU(ED/E_p)<SU(EF)$ and by induction hypothesis
(applied to $D$ and 
$E$), $E$ satisfies $(*)$ over $E_p$. Write
$E=\acl(E'_p,E''_p)$, where $E'_p$ satisfies $(*)$ over $E_p$,
$tp(E'_p/E_p)$ is
hereditarily orthogonal to all types in $\cals(p)$, and
$tp(E''_p/E_p)$ is
$\cals(p)$-analysable. Then $tp(E''_p)$ is
$\cals(p)$-analysable, so that $E''_p=E_p$. On the other hand,
$tp(E/A)$ is $\cals(p)$-analysable, and therefore $E'_p\subseteq
\acl(AE_p)$ (because $tp(E'_p/E_p)$ is hereditarily orthogonal to all
members of $\cals(p)$). Hence $E=\acl(AE_p)$
satisfies $(*)$, a contradiction.

By symmetry,  if $F_p$ is a
 subset of $F$ whose type is $\cals(p)$-analysable, then $F_p\subseteq 
 \acl(\emptyset)$. 

\smallskip\noindent
{\bf Step 4}. $F\subseteq \acl(B_1E)$.

Let $D=\acb(F/B_1E)$. Then $B_1\subseteq D$, and $tp(D/B_1)$ is
almost-$\cals(p)$-internal, because it is contained in the algebraic
closure of $B_1$-conjugates of $F$.
Furthermore, we have $\acb(D/E)=E$: let
$E_0=\acb(D/E)$; from
$E\dnfo_DF$ we deduce $E\dnfo_{DE_0}F$; then $E\dnfo_{E_0}D$ yields
$E\dnfo_{E_0}F$, whence $E_0=E$. We now let $D_1=\acb(E/D)$; then
$D_1\subseteq D\subseteq \acl(B_1E)$,  $tp(D_1/B_1)$ is
almost-$\cals(p)$-internal, and $\acb(D_1/E)=E$ (by Lemma
\ref{lem2}). Since $E$ does not satisfy $(*)$, our 
induction hypothesis implies that either $E\cap D_1\neq \acl(\emptyset)$ or 
$\acl(D_1E)=\acl(EF)$.

Let us  assume that $D_1\cap E\neq \acl(\emptyset)$. Using the
semi-minimal analysis of $tp(D_1\cap E)$, there is $d\in D_1\cap E$ with $tp(d)$
almost-$\cals_i$-internal for some $i$. Since $E_p=F_p=\acl(\emptyset)$, we know
that  $\cals(p)\perp \cals_i$. But $tp(D_1/B_1)$ is 
almost-$\cals(p)$-internal, so that $tp(d/B_1)$
is almost-$\cals(p)$-internal, whence $d\in B_1$. Hence $D_1\cap
E\subset B_1\cap E=\acl(\emptyset)$.

Hence $\acl(D_1E)=\acl(EF)$, which implies $F\subseteq \acl(B_1E)$. 
 
The proof only used the
$\cals(p)$-internality of $tp(F/B_1)$, and 
we  reason in
the same manner with
$\acb(E/AF)$ to
get $E\subseteq \acl(AF)$. Since $E_p=F_p=\acl(\emptyset)$, we know that 
$\cals(p)\perp \cals_1\cup\cdots\cup\cals_m$. 
The final contradiction will come
from the following lemma, taking $\cals=\cals_1\cup\cdots\cup\cals_m$:

\para\vlabel{lem10}{\bf Lemma}. Let $A\subseteq E$ and $B\subseteq F$ be
algebraically closed sets of finite SU-rank such that
$E\cap F=\acl(\emptyset)$, $E$ and $F$ are equi-algebraic over $AB$.
Assume
that for some set $\cals$ of types of
SU-rank $1$, which is closed under $\aut(M/\acl(\emptyset))$-conjugation,
$tp(A)$ and $tp(B)$ are $\cals$-analysable. 
Then $tp(E/\acl(\emptyset))$ and $tp(F/\acl(\emptyset))$ are $\cals$-analysable.

\prf We may assume that $A$ and $B$ are maximal 
subsets of $E$ and $F$ respectively whose type are $\cals$-analysable. If
$E=A$, then $F\subseteq \acl(AB)$, and we are done; similarly if
$F=B$. Assume $E\neq A$, and let $p$ be a type of SU-rank $1$ 
 which is non-orthogonal to $tp(E/A)$; we then let
$E_0\subseteq E$ and $F_0\subseteq F$ be maximal such that $tp(E_0/A)$
and $tp(F_0/B)$ are almost $\cals(p)$-internal. Then $E_0\neq A$ (see
the discussion in \ref{semimin})  and
$p\perp \cals$. If $F_1=\acb(E_0/F)$
and $B_1=\acb(A/F)$, then $tp(F_1/B_1)$ is almost
$\cals(p)$-internal by Lemma \ref{lem1}, so that $F_1\subseteq
F_0$ and $E_0\subseteq \acl(AF_0)$. Similarly, $F_0\subseteq
\acl(BE_0)$. %% Hence we may  assume that $tp(E/A)$ and
%% $tp(F/B)$ are almost-$\cals(p)$-internal for some type $p$ of SU-rank
%% $1$ with 
%% $p\perp \cals$.

We have therefore shown that if the conclusion of the lemma does not
hold, then there is a counterexample $(E,F,A,B)$ where  $tp(E/A)$ and $tp(F/B)$ are almost $\cals(p)$-internal for some 
type $p$ of  
$SU$-rank $1$ which is orthogonal to all members of $\cals$. 
We choose such a counterexample
with $r=SU(B)-SU(B/A)$ minimal. 

Let $E_0=\acb(F/E)$, and $A_0=A\cap E_0$. Then $F\subseteq \acl(BE_0)$,
and $E_0\subseteq \acl(AF)$. Also, $A_0$ is the maximal subset of $E_0$
with an $\cals$-analysable type (by \ref{lem25}), whence $A\dnfo_{A_0}E_0$, and by
transitivity $A\dnfo_{A_0}E_0F$, so that $E_0\subseteq \acl(A_0F)$.
Since $F\neq B$, we have $E_0\neq A_0$, so that $tp(E_0)$ is
not $\cals$-analysable. Replacing $E$ by $E_0$ and $A$ by $A_0$, we
may therefore assume that $E=\acb(F/E)$. (Note that $SU(B/A_0)\geq
SU(B/A)$, so that $SU(B)-SU(B/A_0)\leq SU(B)-SU(B/A)$, and in fact
equality holds by minimality of $r$).

If $r=0$, then $B\dnfo A$; because
$tp(E/A)$ is orthogonal to all types in $\cals$, we obtain  $B\dnfo E$;
since $tp(F/B)$ is 
almost-$\cals(p)$-internal and $E=\acb(F/E)$, we get that 
$tp(E)$ is almost-$\cals(p)$-internal,
a contradiction. Hence $r>0$.

Since $E=\acb(F/E)$, there are $E$-independent realisations  $F_1,\ldots,F_s$
of $tp(F/E)$ such that $E\subseteq \acl(F_1,\ldots,F_s)$. Let $B_i\subset
F_i$ correspond to $B\subset F$. Since
$tp(F_1,\ldots,F_s/B_1,\ldots,B_s)$ is orthogonal to all types in 
$\cals$, we necessarily have $A\subset
\acl(B_1,\ldots,B_s)$. Furthermore, from $B\dnfo_AE$, the sets $B_1,\ldots,B_s$
are independent over $A$. This implies that $\acb(B/A)=A$ by
\ref{sett}(1).

Let $m\leq s$ be minimal such that  $A\subset
\acl(B_1\ldots B_m)$. Then $m>1$ and
$SU(B_m)-SU(B_m/B_1\ldots B_{m-1})<r$ by \ref{sett}(2). We also have $F_m\cap
\acl(F_1\ldots F_{m-1})\subseteq F\cap E=\acl(\emptyset)$, and
$E\subseteq \acl (B_1\ldots
B_mF_i)$ for every $1\leq i\leq m$. Hence $F_1$ and $F_m$ are
equi-algebraic over $\acl(B_1\ldots B_m)$. The induction hypothesis
applied to the quadruple \hfil\break
$(\acl(F_1B_2\dots B_{m-1}), F_m, \acl(B_1\ldots B_{m-1}),B_m)$
gives that $tp(F)$ is
$\cals$-analysable, a contradiction.

\para\vlabel{thm37} Proposition \ref{lem37} has the following
immediate consequence:

\smallskip\noindent
{\bf Theorem}. 
Let $E,F$ be algebraically closed sets, and assume that $SU(E/E\cap F)$
is finite and $F=\acb(E/F)$. Then there are $F_1,\ldots,F_m$
independent over $E\cap F$, types $p_1,\ldots,p_m$  of
$SU$-rank $1$, such that each $tp(F_i/E\cap F)$ is
$\cals(p_i,E\cap F)$-analysable, and $\acl(F_1\ldots F_m)=F$. 

\prf By \ref{sett}, we know that $SU(F/E\cap F)$ is also finite. Replace
$E$ by $E'=\acb(F/E)$; by Lemma \ref{lem2}, $F=\acb(E'/F)$. Then apply
Proposition \ref{lem37} to $E'$, $F$ to get
 the types $p_i$ (which are pairwise orthogonal), and the sets 
$F_i$.

\para\vlabel{rem37}{\bf Remark}. Let $E$, $F$ be as above. Using the semi-minimal
analysis of $tp(F/E\cap F)$, there is some $G=\acl(G)$ independent from
$EF$ over $C=E\cap F$, and a tuple $a\in\acl(GF)$ of realisations of types of
SU-rank $1$ over $G$, such that for any tuple $b$, $b\dfo_GF$ implies
$b\dfo_Ga$ (in other words: $tp(a/G)$ {\em dominates} $tp(F/G)$, see section 5.2 in \cite{[W1]}). Then working over $G$, the
types $p_i$ of Theorem 
\ref{thm37} can be taken to be types over $G$ (see the proof of
\ref{lem37}), and the subsets $F_i$ of $\acl(GF)$ will then realise
$\{p_i\}$-analysable types over $\acl(GC)$. This is slightly stronger
than just saying that the sets $F_i$ realise $\cals(p_i)$-analysable types.

\para\vlabel{prop2-1w} The following result is
similar to Proposition~\ref{prop2-1}. See also Theorem 1.3 in \cite{[MP]}. 

\smallskip\noindent
{\bf Proposition}. Let $\cals$ be a set of types of
SU-rank $1$, which is closed under $\aut(M/\acl(\emptyset))$-conjugation,
let $B$ and $E$ be algebraically closed sets of finite SU-rank, and
assume that $tp(E/B)$ is $\cals$-analysable. Then so is $tp(E/E\cap B)$.

\prf Without loss of generality, $B=\acb(E/B)$. Let $C=E\cap B$, and assume that $tp(E/C)$ is not
$\cals$-analysable. Let $D\subseteq E$ be maximal such that $tp(D/C)$ is
$\cals$-analysable.  As
$B=\acb(E/B)$, Theorem \ref{thm37} gives us two algebraically closed
sets $B_1$ and $B_2$ with $\acl(B_1B_2)=B$, $tp(B_1/C)$
$\cals$-analysable, and $tp(B_2/C)$ $\cals'$-analysable for some
set $\cals'$ of SU-rank $1$ types with algebraically closed base and
such that  $\cals\perp\cals'$. Then
$\acb(D/B)\subseteq B_1$ and $E\dnfo_{D} B_1$ because $tp(E/D)\perp
\cals$ and $tp(B_1/C)$ is 
$\cals$-analysable. If $B_2=C$, then $E\dnfo_DB$, and the
$\cals$-analysability of $tp(E/DB)$ implies the $\cals$-analysability of
$tp(E/D)$, a contradiction. 

Hence $B_2\neq C$, and if $E_2=\acb(B_2/E)$, then $E_2\neq C$ and $E_2$ realises an
$\cals'$-analysable type over $C$. As $E_2\neq C$ and $E\cap B=C$, we have
that $tp(E_2/B)$ is non-algebraic and $\cals'$-analysable. On the other
hand $tp(E_2/B)$ is also $\cals$-analysable because $E_2\subseteq E$,
which gives  
the final contradiction.

\para{\bf Corollary}. Let $\cals$ be a set of types of rank $1$ closed
under $\aut(M/\acl(\emptyset))$-conjugation, and let $E=\acl(E)$ have finite
SU-rank. Then there is $A=\acl(A)\subseteq \acl(E)$ such that $tp(E/A)$ is
$\cals$-analysable, and whenever $B=\acl(B)$ is such that $tp(E/B)$ is
$\cals$-analysable, then $A\subseteq B$.

\prf By \ref{prop2-1w}, it is enough to show that if $A_1, A_2$ are
algebraically closed subsets of $E$ such that $tp(E/A_i)$ is
$\cals$-analysable, then so is $tp(E/A_1\cap A_2)$: but this is obvious,
as $tp(A_1/A_1\cap A_2)$ is $\cals$-analysable, by \ref{prop2-1w}. 

\para{\bf Remark}.  Let $\cals$ be a set of types of rank $1$ closed
under $\aut(M/\acl(\emptyset))$-conjugation, and let $E=\acl(E)$ have finite
SU-rank. Then one can find $\cals'\perp \cals$, closed under $\aut(M/\acl(\emptyset))$-conjugation
and such that $tp(E/\acl(\emptyset))$ is
$(\cals\cup\cals')$-analysable. It follows that for a set $B=\acl(B)$,
$tp(E/B)$ will be hereditarily orthogonal to $\cals'$ if and only if it
is $\cals$-analysable. Thus the above two results can be stated in terms of
hereditary orthogonality to  $\cals'$ instead of $\cals$-analysability.

\medskip
We now state an easy lemma reducing further the problem of showing the
CBP: 

\para \vlabel{lem4}
{\bf Lemma}. Let $\cals$ be a set of types of
SU-rank $1$, which is closed under $\aut(M/\acl(\emptyset))$-conjugation. Assume that there are algebraically closed sets 
$E$ and $F$ whose types over $C=E\cap F$ are $\cals$-analysable,  such
that $\acb(F/E)=E$ and
$tp(E/C)$ is not almost $\cals$-internal. Then there are such sets
$E$ and $F$ whose types over $C$
are $\cals$-analysable in at most two steps, i.e., there is $A\subset
E$
such that $tp(A/C)$ and $tp(E/A)$ are \asi, and similarly for
$F$. Furthermore, $\acb(E/F)=F$. 

\prf We take such a triple $(E,F,C)$ with
    $r=SU(E/C)+SU(F/C)$ minimal, whence $F=\acb(E/F)$.

By the semi-minimal analysis of $tp(F/C)$, there is a proper algebraically
closed subset $B$ of $F$ such that $tp(F/B)$ is \asi. By Lemma
\ref{lem1}, if $A=\acb(B/E)$ then $tp(E/A)$ is \asi. As
$SU(B/C)<SU(F/C)$, the minimality of $r$ implies that $tp(A/C)$ is
\asi. Hence, $A\neq C, E$
because $tp(E/C)$ is not \asi, and $tp(E/C)$ is $\cals$-analysable in
two steps. Since $F$ is contained in
the algebraic closure of realisations of $tp(E/C)$,
$tp(F/C)$
will also be $\cals$-analysable in two steps.

\medskip
We conclude this section with a partial internality result: 

\para\vlabel{lem11}{\bf Lemma}. Let $A\subseteq E$ and $B\subseteq F$ be
algebraically closed sets of finite SU-rank such that
$E\cap F=\acl(\emptyset)$, $E$ and $F$ are equi-algebraic over $AB$.
Assume
that for some set $\cals$ of types of
SU-rank $1$, which is closed under $\aut(M/\acl(\emptyset))$-conjugation,
$tp(A/\acl(\emptyset))$ and $tp(B/\acl(\emptyset))$ are \asi. 
Then $tp(E/\acl(\emptyset))$ and $tp(F/\acl(\emptyset))$ are \asi.

\prf We work over $\acl(\emptyset)$. By Lemma \ref{lem10}, we already know that $tp(E)$ and $tp(F)$ are
$\cals$-analysable. Hence, reasoning as in the first paragraph of the
proof of \ref{lem10}, we may assume that $A$, $B$ are maximal subsets of
$E$ and $F$ respectively which realise \asi-types, and that $tp(E/A)$
and $tp(F/B)$ are \asi, but neither $tp(E)$ nor $tp(F)$ is \asi. The maximality of $A$ and $B$ implies that
$A\dnfo_BF$ and $E\dnfo_AB$. 

Working over some $C=\acl(C)$, independent from $EF$, and using
\ref{sett}(3) and \ref{lem25}, we may assume that $A$ is the algebraic
closure of a tuple of realisations of types in $\cals$. We choose a
counterexample $(E,F,A,B)$ with $r=SU(A)-SU(A/B)$ minimal. 
 If $r=0$, then $A\dnfo B$ so that
$A\dnfo F$. Letting $F_0=\acb(E/F)$,  this implies that $F_0$ realises an
\asi-type, hence is contained in $B$. But $E\dnfo_{F_0}F$ then implies
$F\subset B$, which is absurd. So we may assume that $r>0$.

Let $(F_iB_i)_{i>0}$ be a sequence of $E$-independent realisations of
$tp(FB/E)$. If $A_0=\acb(B/A)$, then for some $s>1$, we
have $\acl(B_1\ldots B_s)\supset A_0$, and we take a minimal such
$s$. As $A$ is the algebraic closure of realisations of types of SU-rank
$1$, there is a finite tuple $a\subset A$ such that $a\dnfo A_0$ and
$\acl(A_0a)=A$. Then $a\dnfo A_0B$, which implies $a\dnfo A_0F$ (by
transitivity and because $A\dnfo_BF$). 

Let $A'=\acl(aB_1)$, $B'=\acl(B_2\ldots B_s)$, $E'=\acl(A'F_1)$,
$F'=\acl(B'F_s)$. Then $a\dnfo A_0F$, and 
$a\dnfo B_1F'$. Hence in particular, $A'\cap F'=\acl(\emptyset)$ (use
\ref{sett}(3) and the fact that $B_1\cap B'=\acl(\emptyset)$). 
Moreover, $SU(A')-SU(A'/B')=SU(B_1)-SU(B_1/B_2\ldots
B_s)<r$ by \ref{sett}(2), and $E'$ and $F'$ are equi-algebraic over $A'B'$.  
In order to reach a contradiction, it therefore suffices to show that
$E'\cap F'=\acl(\emptyset)$: our induction hypothesis gives that
$tp(E/\acl(\emptyset)$ is \asi, which  implies that
$tp(F_1/\acl(\emptyset))=tp(F/\acl(\emptyset))$ is also \asi.

By Lemma \ref{lem25}, $A'$ and $B'$ are maximal
subsets of $E'$, $F'$ respectively which realise \asi-types.  
Assume $E'\cap F'\neq \acl(\emptyset)$; by the semi-minimal analysis
of $tp(E'\cap F'/\acl(\emptyset))$,
there is $d\in E'\cap F'$ realising
an \asi-type. Then $d\in B'\cap A'\subseteq F'\cap A'=\acl(\emptyset)$, which
gives us the desired contradiction.

\sect{Further properties of theories with the CBP}

{\bf Description of the  results of this section}. Assumptions on
$M$ and $T$ are as in the previous section: $T$ is supersimple and
eliminates imaginaries, $M$ is sufficiently saturated.  Most results are
proved under the additional hypothesis of the CBP. We start by
proving one of the main results of the paper:  

\smallskip\noindent
{\bf Theorem \ref{prop2}} (CBP). If  $E$ and $F$
are algebraically closed sets of finite SU-rank 
over their intersection $C$ and are such that $E=\acb(F/E)$, then
$tp(E/C)$ is \asi{} for
some family $\cals$ of types of SU-rank $1$.

\smallskip
Note that under the same hypotheses, if $tp(E/F)$
is $\cals'$-analysable for some set $\cals'$ of types of SU-rank $1$
with algebraically closed base, then 
$tp(E/C)$ is almost-$\cals'$-internal. We then show 

\smallskip\noindent
{\bf Theorem \ref{thm2-1}} (CBP). Assume that $E=\acl(E)$ has finite
SU-rank, and let $\cals$ be 
 a collection of types of SU-rank 1,
closed under conjugation. Then there is $A=\acl(A)\subseteq E$ such that $tp(E/A)$
is \asi, and whenever $B=\acl(B)$ is such that $tp(E/B)$ is \asi, then
$B\supseteq A$. 

An immediate consequence of \ref{thm2-1} is that the CBP implies the
Uniform CBP (UCBP); this answers a question of Moosa and Pillay
\cite{[MP]}. We end the section with three results, which
were proved with an eye towards geometric applications. The first result
is valid in a general setting (as will be clear from the proof), and can
be viewed as showing the existence of a ``largest internal
quotient''; the second can be viewed as showing the existence of a
``maximal internal fiber'', and the third one as a descent result.

\para\vlabel{prop2}{\bf Theorem} (CBP). 
If   $E$ and $F$ are algebraically closed sets of 
finite SU-rank over their
intersection $C$ 
and $E=\acb(F/E)$, then $tp(E/C)$ is \asi\ for some family $\cals$ of types
of SU-rank $1$.

\prf We assume the result false. By Lemma 
\ref{lem4}, 
there is a 
counterexample $(E,F,C)$ with  $tp(E/C)$, $tp(F/C)$
$\cals$-analysable in two steps, and which also satisfies $\acb(E/F)=F$.
By Theorem  \ref{thm37} (see also \ref{rem37} and 
\ref{sett}(3)), working over a larger set $G=\acl(G)$, we can write $E$
as $\acl(E_1\ldots
E_m)$ for some sets $E_i$  which are independent over $G$, realise
$\{p_i\}$-analysable 
types over $G$, and some $E_i$ will not realise an
almost-$\{p_i\}$-internal type over $G$. 

Hence, we may assume that 
$\cals=\{p\}$ for some type $p$ of SU-rank $1$. {\em  For ease of
  notation we will   assume that the language
  contains constant symbols for the elements of $G$.} 

Let $A_0\subset E$ and $B\subset F$ be maximal realizing \asi\ types, so that
$tp(E/A_0)$ and $tp(F/B)$ are \asi. Then $E\dnfo_{A_0}B$ since
$\acb(B/E)$ is \asi\ and therefore contained in $A_0$, and similarly
$A_0\dnfo_BF$. Enlarging $G$ (and using \ref{sett}(3)), we may 
assume that $A_0$ and $B$ are the algebraic closures of tuples of
realisations of 
$p$.
    Let $A=\acb(B/E)$. Then  $A\subseteq A_0$, and
    $tp(E/A)$ is \asi{}
(by \ref{lem1}).

The proof is by induction on $r=SU(B)-SU(B/A_0)$
($=SU(A_0)-SU(A_0/B)$). If $r=0$, then $A_0\dnfo B$ and from
$tp(F/B)$ \asi{} and $\acb(F/E)=E$ we deduce that $tp(E)$ is \asi, a
contradiction. Hence $r>0$.

\smallskip\noindent
{\bf Step 1}. We may assume $A_0=A$.

We know that $A_0=\acl(a_0)$ for some
tuple $a_0$ of realisations of $p$; take $a\subseteq a_0$ maximal
independent over $A$ and such that $a\dnfo A$. Then  $$\acl(Aa)=A_0, \qquad
 a\dnfo AF \qquad\hbox{and}\qquad
SU(B/a)-SU(B/A_0)=r.$$ Furthermore, $tp(E/a)$ is not \asi:
otherwise, $\acb(E/Fa)$ would also be \asi, hence contained in $\acl(Ba)$
by maximality of $B$ (see \ref{lem25}(1)); from
$E\dnfo_{Ba}F$ and $F\dnfo_BA_0$ we would then deduce $E\dnfo_BF$,
i.e. $F=B$, which is absurd. 
We 
will now show that $E\cap \acl(Fa)=\acl(a)$. Enlarging
$G$ this will allow us to assume  $A=A_0$.

Let $D=E\cap \acl(Fa)$. Since $a\dnfo AF$, $A_0\cap \acl(Fa)=\acl(a)$ by
\ref{sett}(3). 
The set $\acb(B/D)$ is \asi, hence contained in $A_0\cap \acl(Fa)=\acl(a)$,
so that $B\dnfo D$ because $F\dnfo a$. From $D\subset \acl(Fa)$ and
the
\asi ity of $tp(Fa/B)$ we obtain that $tp(D/B)$ is \asi, and therefore
also  $tp(D)$, so that
$D\subseteq A_0\cap \acl(Fa)=\acl(a)$.

\smallskip\noindent
{\bf Step 2}. We may assume  $E\subseteq \acl(AF)$.

By assumption, there is an algebraically closed set $J$ containing $F$,
such that $J\dnfo_FE$,  and a tuple $g$ of
realisations of $p$ such that $E\subseteq \acl(Jg)$. 
Then there is a subset $e$ of $g$, consisting of independent tuples over
$AJ$, and such that $E\subseteq \acl(AJe)$ and $e\dnfo AJ$.  Since $J\supseteq F$ and
$e\dnfo J$, we then
have $\acb(Je/Ee)=\acl(Ee)$ and $\acb(Be/Ae)=\acb(Be/Ee)=\acl(Ae)$ (use
$\acb(J/E)=E$ and \ref{sett}(1)).

\smallskip\noindent
{\bf Claim}. $\acl(Ee)\cap \acl(Je)=\acl(e)$.

Let $D=\acl(Ee)\cap \acl(Je)$. Since $e\dnfo AJ$, we obtain 
$\acl(Ae)\cap  \acl(Je)=\acl(e)$ by \ref{sett}(3).

We know that if $D_0=\acb(A/D)$, then $tp(D_0)$ is \asi; the maximal
\asi{} subset of 
$\acl(Ee)$ is $\acl(Ae)$ by \ref{lem25}(1), and therefore $D_0\subseteq \acl(Ae)\cap
\acl(Je)=\acl(e)$. Hence $A\dnfo D$, and $tp(D)$ is \asi{}  (because
$D\subset \acl(Ee)$
and $tp(Ee/A)$ is \asi). Reasoning as we did for $D_0$, we obtain
$D\subseteq \acl(Ae)\cap 
\acl(Je)=\acl(e)$.

  From $E\dnfo_FJ$, $e\dnfo AJ$ and $A\dnfo_BF$ we deduce
$$SU(A/e)=SU(A)\hbox{ and }
SU(A/Je)=SU(A/J)=SU(A/F)=SU(A/B),$$ so that
$$SU(A/e)-SU(A/Je)=SU(A)-SU(A/B)=SU(B)-SU(B/A)=r.$$Because $e\dnfo A$ and by 
maximality of $A$, we get $e\dnfo E$; thus $tp(E/e)$ is not \asi. As we
saw above, we have $\acb(Je/Ee)=\acl(Ee)$. 
Hence,  working over $\acl(e)$ and replacing $F$
by $J$, we may assume $E\subseteq \acl(AF)$. 

\smallskip\noindent
{\bf Step 3}. The final contradiction.

Let $(F_nB_n)_{n\in\nat}$ be a sequence of $E$-independent realisations
of $tp(FB/E)$. From $B\dnfo_AE$, it follows that the sets $B_n$ are
independent over $A$. By \ref{sett}(1), and because
$A=\acb(B/A)$, there is $m$ such that $A\subset \acl(B_1\ldots B_m)$;
take the minimal such $m$. Then $E\subset \acl(B_1\ldots B_mF_i)$ for every
$i$, so that in particular $F_1\dfo_{B_1\ldots B_m}F_m$. On the other
hand, we know that $F_1\cap \acl(B_2\ldots
B_{m-1}F_m)\subseteq F\cap E=\acl(\emptyset)$, and
$SU(B_1)-SU(B_1/B_2\ldots B_m)<r$ by
minimality of $m$. We apply the induction hypothesis to $(F_1,B_1)$ and
$(\acl(B_2\ldots B_{m-1}F_m), \acl(B_2\ldots B_m))$: if $J=\acb(B_2\ldots
B_{m-1}F_m/F_1)$, then $J\not\subseteq B_1$ and $tp(J)$ is \asi.  This
contradicts the maximality of $B$, and finishes the proof.

\para We will now prove some more results for supersimple theories with
the CBP. Note that Proposition \ref{prop2-1} below implies Theorem
\ref{prop2} and is therefore equivalent to it. It was first proved by
Moosa and Pillay in the stable context,  see \cite{[MP]}.%, but we will give a shorter proof here.

\para\vlabel{prop2-1}{\bf Proposition} (CBP). Let $B$ and  $E$ be algebraically closed sets,
with $SU(E)<\infty$, and assume that $tp(E/B)$ is \asi, for some
collection $\cals$ of types of SU-rank $1$, which is closed under
$\aut(M/\acl(\emptyset))$-conjugation. Then  $tp(E/E\cap B)$ is \asi.

\prf Let $C=B\cap E$, and let $A\subseteq E$ be maximal such that
$tp(A/C)$ is \asi. If $B_0=\acb(E/B)$, then $tp(E/B_0)$ is also \asi,
and we may therefore assume that  $B=\acb(E/B)$. By Proposition
\ref{prop2-1w}, we know that $tp(E/C)$ is $\cals$-analysable, and this
implies that $tp(B/C)$ is also $\cals$-analysable. On the other hand, by 
Theorem
\ref{prop2}, $tp(B/C)$ 
is almost-$\cals'$-internal, for some collection $\cals'$ of types of
SU-rank $1$ containing $\cals$, and these two facts imply that $tp(B/C)$
is \asi. 

Assume $E\neq A$. By assumption, there is some
$F=\acl(F)$, independent from $E$ over $B$, and such that  $E$
is equi-algebraic over $F$ with some finite tuple of realizations of types in
$\cals$. 

\smallskip\noindent
{\bf Claim}. $\acl(AF)\neq \acl(EF)$.

Otherwise, $A\subseteq E$ and $E\dnfo_BF$ would imply  $E\subseteq
\acl(AB)$. As $tp(B/C)$ 
is almost-$\cals$-internal, this would imply that also $tp(E/C)$ is
almost-$\cals$-internal, a contradiction. 
 
We may therefore choose some $e\in\acl(EF)\setminus \acl(AF)$ which
realises a type in $\cals$. Then
$E_0=\acb(Fe/E)\not\subseteq A$, since $e\in \acl(FE_0)\setminus
\acl(FA)$. Note that $E\cap F=E\cap B=C$. 

Let $D=\acl(Fe)\cap E$. Then $D\cap F=C$, and by Theorem \ref{prop2}
$tp(E_0/D)$ is \asi\ (because $E_0\subseteq E$ and $tp(E/C)$ is $\cals$-analysable). If $D=C$, this gives
us the desired contradiction, as $E_0\not\subseteq A$, and $A$ was
maximal contained in $E$ with $tp(A/C)$ \asi.

Assume therefore that $D\neq C$. Then $SU(D/F)=1$, because $SU(e/F)=1$
and $D\subset \acl(Fe)$. If $D\dnfo_CF$, then
$SU(D/C)=1$, which implies that $tp(D/C)$ is \asi. In that case we let
$D_0=D$. If $D\dfo_CF$, we define $D_0=\acb(F/D)$; then $tp(D_0/D\cap
F)$ is \asi\ by \ref{prop2}. Hence, as $D_0\subseteq \acl(Fe)$, and $D_0\not\subseteq F$,
we have that $e\in\acl(FD_0)$, and $tp(D_0/C)$ is \asi. As $D_0\subseteq
D\subseteq E$, we obtain $D_0\subseteq A$, whence $e\in \acl(FA)$, which gives us the desired
contradiction and finishes the proof.

\para\vlabel{lem2-1}{\bf Lemma}  (CBP). Let $E=\acl(E)$ be of finite SU-rank over some $C=\acl(C)$, and let $\cals$ be a collection of types of SU-rank 1,
closed under $\aut(M/\acl(\emptyset))$-conjugation. Assume that $A_i=\acl(A_i)\subset E$, $i=1,2$, are such
that $A_1\cap A_2=C$, and $tp(E/A_i)$ is \asi\ for $i=1,2$. Then
$tp(E/C)$ is \asi. 

\prf Let $A\subseteq E$ be maximal such that $tp(A/C)$ is \asi. Then
$A_1\subseteq A$: by hypothesis, $tp(A_1/A_2)$ is \asi, and by
\ref{prop2-1}, $tp(A_1/C)$ is \asi. Reasoning similarly with $A_2$, we
obtain that $A_1A_2\subseteq A$. If $F=\acl(F)\supset C$ is independent
from $E$ over $C$, and 
$tp(E/F)$ is \asi, then so is $tp(E/C)$, and we may therefore extend
$C$, to a larger set over which $A$ is equi-algebraic with a tuple of realisations
of types in $\cals$ (by Lemma \ref{lem25}(2), we will not lose the
maximality of $A$). Hence, we may assume that in $A$ there is a tuple $a$ of
realisations of types in $\cals$ such that $a\dnfo_CA_1A_2$,
and $A=\acl(CA_1A_2a)$. Note that we still have
$\acl(A_1a)\cap A_2=C$: since $a\dnfo_CA_1A_2$, we know by (1.1)(3) %\ref{sett}
that $\acl(A_1a)\cap \acl(A_2a)=\acl(Ca)$; hence $\acl(A_1a)\cap
A_2\subseteq \acl(Ca)\cap A_2=C$. Thus, 
replacing $A_1$ by $\acl(A_1a)$
we may assume that $\acl(A_1A_2)=A$. 

By assumption, for $i=1,2$, there are $F_i=\acl(F_i)$ containing $A_i$, independent from $E$
over $A_i$, and such that  $E$ is equi-algebraic  over $F_i$ with some
tuple $b_i$ of realisations of types in $\cals$. We may choose  $F_2$
independent from $EF_1$ over $A_2$; then $F_1\dnfo_E F_2$, whence also $F_1$ is independent from
$EF_2$ over $A_1$, and  
$$C=A_1\cap A_2=F_1\cap F_2; \qquad \acl(F_1b_1)\cap F_2=A_2; \qquad
F_1\cap \acl(F_2b_2)=A_1$$ 
(use $\acl(F_ib_i)=\acl(F_iE)$, $E\cap F_j=A_j$).
For $i=1,2$, choose $e_i\subset b_i$ maximal
independent over $F_iA$. Then $E\subseteq \acl(F_iAe_i)$, and $A\cap \acl(F_ie_i)=A_i$. Furthermore
$$\acl(F_1e_1)\cap F_2=F_1\cap\acl(F_2e_2)=A_1\cap A_2=C.$$
Let $D_0=\acl(F_1e_1)\cap \acl(F_2e_2)$. As $D_0\subseteq \acl(F_1e_1)$,
$tp(D_0/F_1)$ is \asi;
by \ref{prop2-1}, $tp(D_0/D_0\cap F_1)$ is \asi; hence $tp(D_0/C)$ is
\asi\ because
$\acl(F_2e_2)\cap F_1=C$, and this implies that
$D_0\cap E\subseteq A$. Therefore 
$$D_0\cap E=D_0\cap A=\acl(F_1e_1)\cap \acl(F_2e_2)\cap A=A_1\cap A_2=C.$$ 
Let $D_1=\acb(F_1e_1/F_2e_2)$. Then $tp(D_1/D_0)$ is \asi, by \ref{prop2}. We know that
$F_1e_1$ and $F_2e_2$ are independent over $D_1$, and therefore 
$$F_1e_1\dnfo_{D_1A_1A_2}F_2e_2.$$
Since $\acl(A_1A_2)=A$ and $E\subseteq \acl(F_iAe_i)$, we get
$E\subseteq \acl(D_1A)$. Hence $tp(E/D_0)$ is \asi, and so is
$tp(E/D_0\cap E)$ (by \ref{prop2-1}). As $D_0\cap E=C$, we get the result. 

\para{\bf Theorem}\vlabel{thm2-1} (CBP). Assume that $E=\acl(E)$ has finite
SU-rank, and let $\cals$ be 
 a collection of types of SU-rank 1,
closed under $\aut(M/\acl(\emptyset))$-conjugation. Then there is $A=\acl(A)\subseteq E$ such that
 $tp(E/A)$ 
is \asi, and whenever $B=\acl(B)$ is such that $tp(E/B)$ is \asi, then
$B\supseteq A$.

\prf This follows immediately from Proposition \ref{prop2-1} and
 Lemma \ref{lem2-1}. 

\para{\bf Theorem}\vlabel{thm2-2} (CBP). Let $B=\acb(A/B)$, where
$A=\acl(A)$ has finite SU-rank, and let
$\cals$ be  a collection of types of SU-rank 1,
closed under $\aut(M/\acl(\emptyset))$-conjugation, and such that $tp(B/A)$ is \asi. If
$C=\acl(C)$ is such that $tp(A/C)$ is \asi, then so is $tp(AB/C)$. That
is, $T$ has the UCBP. 

\prf Let $D=\acb(B/A)$. Then $tp(D/B)$ is \asi, and $B=\acb(D/B)$
(by Lemma \ref{lem2}). 
As $D\subseteq A$, $tp(D/C)$ is \asi, and by \ref{lem2-1}, so is $tp(D/B\cap
C)$; 
this implies that 
$tp(B/C)$ is also \asi, since $B$ is contained in the algebraic
closure of realisations of the \asi-type $tp(D/B\cap C)$ (see \ref{sett}(1)). 

\para\vlabel{gp1}{\bf Proposition} (CBP). Let $G$ be a group of finite SU-rank,  let
$p$ be a type (over $\emptyset$) realised by $a\in G$, and let
$H=\Stab(p)$ be the left stabilizer of $p$. If $d$ is the code of
$H\cdot a$, then $tp(d)$ is \asi, where
$\cals$ is the collection of non-locally modular types of SU-rank $1$
and with algebraically closed base. 

\prf The proof is essentially identical to the one given in \cite{[PZ]},
Corollary 
3.11, where it was done in the stable case. 
Let
$c\in G$ be a generic of $G$ over $a$, and let $D=\acb(c/e)$, where
$e=a\cdot c$. 

We will first show  that $d\in\acl(Dc)$. By genericity of $c$, we
know that $e\dnfo a$, and therefore $a\dnfo De$. The set $D$ has the following property: if
$e_1,e_2$ are $D$-independent realisations of $tp(e/D)$, then there is
$c'$ independent from $e_1e_2$ over $D$, and such that
$tp(c'e_i/D)=tp(ce/D)$ for $i=1,2$. If $a_i=e_i\cdot{c'}\inv$, then
$e_1\cdot e_2\inv
=a_1\cdot a_2\inv$, and $a_1,a_2$ realise $p$. We then deduce successively the
following relations: \\
\begin{gather*}
c'\dnfo_D{e_1\,e_2};\ \ \ c'\dnfo_D (e_1\cdot e_2\inv)\, e_2; \ \ \ c'\dnfo_{De_2}
e_1\cdot e_2\inv;\ \ \ a_2\dnfo_{De_2}e_1\cdot e_2\inv;\\
\hbox{since }a_2\dnfo D e_2,\ \hbox{transitivity implies
}a_2\dnfo D e_1\cdot e_2\inv.
\end{gather*}
As both $a_1$ and $a_2$ realise $p$, and
$a_1=e_1\cdot e_2\inv \cdot a_2$, we get that $e_1\cdot e_2\inv\in H$. 
So we have shown that if $e_1$, $e_2$ are any $D$-independent
realisations of $tp(e/D)$, then $e_1\cdot e_2\inv\in H$. Hence, if $e_1$
and $e_2$
realise $tp(e/D)$, then $e_1\cdot e_2\inv \in H$.

If $\tau\in\aut(M/Dc)$, then
$\tau(e)\cdot e\inv \in H$, and $\tau(a)=\tau(e)\cdot e\inv \cdot a\in
H\cdot a$. This shows
that $d\in\acl(Dc)$.

By the CBP, we know that $tp(D/\acl(c))$ is \asi, and therefore so is
$tp(d/\acl(c))$. But on the other hand, we know that $d\in\acl(a)$ and
$a\dnfo c$: hence $d\dnfo c$ and $tp(d)$ is \asi.

\para\vlabel{gp2}{\bf Remark/Corollary} (CBP). Let $G$ be a group of finite SU-rank,
and $p$ a 
type over $\emptyset$, realised by $a\in G$. Let $b\in \dcl(a)$ be
maximal realising an \asi-type, and let $S=\{g\in G\mid
tp(g \cdot a/b)=tp(a/b)\}$, and let $N$ be the subgroup of $G$ generated by
$S$. Then $N\subseteq H$, where $H$ is the left stabiliser of $p$.

\prf If $\pi:G\to H\backslash G$ is the natural projection, then we know
that $H\cdot a$ is coded by $\pi(a)$. By \ref{gp1}, $tp(\pi(a))$ is \asi, and
therefore $\pi(a)\in\dcl(b)$. By definition of $b$, $tp(a'/b)=tp(a/b)$
implies $a'\in H\cdot a$ and $a'\cdot a\inv\in S$, which gives the result.

\para The next results allow us in many cases to pass from the
algebraic closure of a set to the set itself. In
geometric situations, it will allow us to replace correspondences by
rational maps.  The delicate point is that in general, if 
$B=\acl(B)\subset \acl(A)$ and $B_0=B\cap A$, it may happen that
$B\neq \acl(B_0)$. 
The first result, Observation \ref{intquot}, does not
need the CBP hypothesis. 

In what follows, we work over $\emptyset$, and have a set $\cals$
of SU-rank $1$ 
types with algebraically closed base, and
which is closed under $\aut(M)$-conjugation.

\para\vlabel{intquot} {\bf Observation}. Let $a$ a tuple, let $B=\acl(B)$ be maximal
contained in $\acl(a)$ and such that $tp(B)$ is \asi. Let
$B_0=\dcl(a)\cap B$; then $\acl(B_0)=B$. 

\prf Let $b\in B$ be such that $B=\acl(B)$, and let $b'$ be a conjugate
of $b$ over $\dcl(a)$. Then $tp(b')=tp(b)$, and therefore $tp(b')$
is \asi. Hence, if $c$ is a  tuple
encoding the set of conjugates of $b$ over $\dcl(a)$, then
$c\in \dcl(a)$, and $tp(c)$ is \asi, so that $c\in B_0$. As
$b\in \acl(c)$, we get $\acl(c)=B$.

\para\vlabel{intfiber} {\bf Proposition} (CBP). Let $a$ be a tuple of
finite SU-rank, let
$B=\acl(B)$ be  
such that $tp(a/B)$ is \asi. If $B_0=B\cap \dcl(a)$, then $tp(a/B_0)$
is \asi.

\prf We may assume that $B$ is minimal algebraically closed  such that
$tp(a/B)$ is \asi. Choose a tuple $b\in B$ such that $B=\acl(b)$. 
If $b'$ is a conjugate of
$b$ over $\dcl(a)$, then
$tp(a,b)=tp(a,b')$, and therefore $tp(a/b')$ is also \asi. The
minimality of $B$ (and \ref{lem2-1}) implies that
$\acl(b')=\acl(b)$. Hence,  if $c$
is a tuple encoding the set of conjugates of $b$ over
$\dcl(a)$, then $\acl(c)=\acl(b)$; as 
$c\in B\cap \dcl(a)=B_0$, we get
 $B=\acl(B_0)$.

\para\vlabel{descent}{\bf Proposition} (CBP). Let $a_1,a_2,b_1,b_2$ be tuples
of finite SU-rank and assume that
\begin{itemize}
\item[$\bullet$]{$tp(b_2)$ is \asi,}
% \item{}{$tp(a_i/b_i)\vdash tp(a_i/\acl(b_i)$ for $i=1,2$ and has finite
% rank,}
\item[$\bullet$]{$\acl(b_1)\cap \acl(b_2)=\acl(\emptyset)$,}
\item[$\bullet$]{$a_1\dnfo_{b_1}b_2$  and $a_2\dnfo_{b_2}b_1$, }
\item[$\bullet$]{$a_2\in \acl(a_1b_1b_2)$.}
\end{itemize}
Then there is $e\subset \dcl(a_2b_2)$ such that $tp(a_2/e)$ is \asi\ and
$e\dnfo{b_2}$. In particular, if $tp(a_2/b_2)$ is hereditarily
orthogonal to all types in $\cals$, then $a_2\in \acl(eb_2)$. 

\prf If
$C=\acb(a_1b_1/a_2b_2)$, then $a_2\in \acl(Cb_2)$. Let
$D=\acl(a_1b_1)\cap \acl(a_2b_2)$. As $D\subset \acl(a_ib_i)$ for $i=1,2$, we have
$D\dnfo_{b_1}b_2$ and $D\dnfo_{b_2}b_1$. Hence $D\dnfo b_1b_2$ because
$\acl(b_1)\cap \acl(b_2)=\acl(\emptyset)$. Furthermore, we know by
Theorem \ref{thm37} 
that 
there is a set $\cals'$ of SU-rank $1$ types orthogonal to all members
of $\cals$ and such that 
$tp(C/D)$ is almost-$(\cals\cup\cals')$-internal. We may write $C$ as
$\acl(c_1c_2)$ where $tp(c_1/D)$ is \asi, and $tp(c_2/D)$ is
almost-$\cals'$-internal. Then $\acl(c_2D)\dnfo b_2$ because
$tp(c_2/D)$ is hereditarily orthogonal to all members of $\cals$ and
$tp(b_2)$ is \asi. Furthermore, as $a_2\in\acl(Dc_1c_2b_2)$, it follows
that $tp(a_2/\acl(Dc_2))$ is \asi. Now, 
$Dc_2\subseteq \acl(a_2b_2)$, and Proposition \ref{intfiber} implies that if
$e=\acl(Dc_2)\cap \dcl(a_2b_2)$, then $tp(a_2/e)$ is \asi. 

The last assertion is clear: $tp(a_2/e)$ \asi\ implies
$tp(a_2/eb_2)$ \asi, and our assumption of hereditary orthogonality implies that
$a_2\in \acl(eb_2)$.

\para{\bf Concluding remarks}. Inspection of the proofs shows that our
assumption of supersimplicity on the ambient theory is unnecessary, as
long as one restricts 
one's attention to types ranked by the SU-rank, and the relevant 
 hyperimaginaries and imaginaries are eliminated. Thus, the results of
 section 1 do apply to types of finite U-rank in separably closed fields
 of finite degree of imperfection. It is unknown whether this family of
 types enjoys
 the CBP, we will explain now  what one needs to prove. Let $K$ be a
 separably closed field of characteristic $p>0$ and finite (positive)
 degree of imperfection. It follows from results of Margit Messmer, 
Hrushovski and Fran\c coise Delon (see e.g. \cite{B}), that  a type of finite
U-rank which is not one-based is non-orthogonal to the generic
type $q$ of $\bigcap_nK^{p^n}$. By \ref{thm37}, it is therefore enough to show
that the family of all $\{q\}$-analysable types has the CBP. A partial
result in this direction is obtained by 
Pillay and Ziegler  in \cite{[PZ]}: they show that the family of
 very thin types has the
 CBP. Thus, the results of section 2 apply for the family $\calp$ of
 very thin types. Unfortunately, Pillay and Ziegler also 
  give an example of a $\{q\}$-analysable type (of U-rank $2$) which is
  not very thin.

\smallskip
The result of Pillay and Ziegler on types in differentially closed
fields of characteristic $0$ is stronger than the CBP: indeed, if $\Cb$
denotes the usual canonical base, then they show that given two tuples
$a$ and $b$
of finite rank such that $b=\Cb(a/b)$, then $tp(b/a)$ is internal to the
constants. It would be interesting to know  whether this implies
that $tp(b/C)$ is also internal to the constants (as opposed to almost-internal to the constants), under some reasonable conditions
on $a$, $b$, and with $C=\acl(a)\cap \acl(b)$, or even  $C=\dcl(a)\cap
\dcl(b)$.

%\eject
\sect{Existentially closed difference fields have the CBP}

\medskip
Recall that a difference field is a field with a distinguished
endomorphism (usually denoted by $\si$),  which we study in the
language of rings augmented by a symbol for $\si$. A difference field
$K$ is {\em inversive} if $\si(K)=K$. We refer to \cite{[C]} for basic
algebraic results on difference fields, and to \cite{[CH]} for basic
model-theoretic results. 
 Any completion of the theory ACFA of existentially closed  difference
fields
is supersimple and eliminates imaginaries. Moreover, if $K$ is an
existentially closed 
difference field, and $A\subseteq K$, then $\acl(A)$ is the smallest
algebraically closed  subfield $B$ of $K$ satisfying
$\si(B)=B$ and containing $A$.
Independence of algebraically closed sets coincides with independence
in the
sense of the theory of algebraically closed fields, i.e., if
$C\subseteq
A,B$ are algebraically closed difference subfields of $K$, then $A$ and $B$ are
independent over $C$ if and only if $A$ and $B$ are linearly disjoint
over
$C$.

As an immediate corollary of the results of Pillay and Ziegler and of
Proposition \ref{prop2}, we then obtain

\para\vlabel{prop3}{\bf Proposition}. Let $(K,D)$ [resp.  $(K,\si)$]
be
a differentially closed field [resp., an existentially closed  difference field] of
characteristic $0$.  Let $C\subseteq A,B$ be algebraically closed
differential [resp. difference] subfields of $K$, with
$SU(B/C)<\omega$. Assume that $A=\acb(B/A)$. Then $tp(A/A\cap B)$
is almost internal to $Dx=0$ [resp.  $\si(x)=x$].

\para
{\bf Notation}.
We denote by $A^{alg}$ the field-theoretic algebraic closure of a
field
$A$, and by $A^s$ its separable closure. If $\si(E)=E$ is a
difference subfield of
the inversive difference field $K$, and $a$ is a tuple of elements of $K$, then
$E(a)_\si$ denotes the (inversive) difference subfield $E(\si^i(a)\mid i\in\zee)$
of $K$.
If $\tau$ is an automorphism of $K$, we denote by $\Fix(\tau)$ the
subfield of $K$ consisting of elements fixed by $\tau$. We denote by
$\Frob$ the Frobenius map $x\mapsto x^p$. 

\para{\bf $p$-bases and degree of imperfection}. For details and
proofs, see \cite{[B]} \S 13. Let $E\subseteq K\subseteq L\subseteq K^{alg}$ be fields
of characteristic $p>0$, 
with $E$ perfect and $tr.deg(K/E)=d<\infty$. Then $[K:K^p]=p^e$ for
some $e\leq d$, 
and there is an $e$-tuple $c$ of elements of $K$ such 
that $K=K^p[c]$. Such a tuple is called a {\em $p$-basis of $K$} and
its elements are algebraically independent over $E$.  Moreover, if
$e=d$, then $c$ is a {\em separating transcendence basis} of $K$ over $E$,
i.e., $K\subseteq E(c)^s$. The integer $e$ is
called the {\em degree of imperfection of $K$}.

We also have: $[L:L^p]$ divides $p^e$, and $[L:L^p]=p^e$  if
$L\subseteq K^s$ or if $[L:K]<\infty$.

\para\vlabel{lem38}{\bf Lemma}. Let $(K,\si)$ be an existentially closed  difference
field of
characteristic $p>0$, let $E=\acl(E)\subset K$, $a$ a finite tuple in
$K$, and assume that
$tp(a/E)$ is $\Fix(\si)$-analysable. Then there is a finite tuple $b$
such that
$E(a)_\si=E(b)_\si$, and
$\si(b),\si\inv(b)\in E(b)^s$.

\prf We will  show that if $d=tr.deg(E(a)_\si/E)$, then
$[E(a)_\si:E(a^p)_\si]=p^d$ and $d<\infty$. This will yield the result: let $c$ be a
$p$-basis of $E(a)_\si$. Then $E(a)_\si \subseteq E(c)^s$, and
therefore $E(a)_\si=E(c,a)_\si\subseteq E(c,a)^s$.

The proof is by induction on the length of a semi-minimal analysis in
$\Fix(\si)$ of
$tp(a/E)$. Assume first that $tp(a/E)$ is almost-$\Fix(\si)$-internal.
Let $F=\acl(F)$ be independent from $a$ over $E$, and such that $a$
is equi-algebraic over $F$ with some finite tuple $b$ of $\Fix(\si)$. We may
assume
that  $a\in F(b)^{s}$ (we replace $b$ by $b^{1/p^n}$ if
necessary). From $\si(b)=b$, we deduce that $F(a)_\si\subseteq F(b)^s$,
and
therefore $$p^d\geq [F(a)_\si:F(a^p)_\si]\geq [F(b):F(b^p)]=p^d.$$ As $F$ was
linearly
disjoint  from $E(a)_\si$ over $E$, this shows
$[E(a)_\si:E(a^p)_\si]=p^d$, with $d<\infty$.

For the general case, choose  $a_1,\ldots,a_n\in \acl(Ea)$ such that
$a\in E(a_1,\ldots,a_n)_\si$, and
for every $i$, $tp(a_i/\acl(Ea_1,\ldots,a_{i-1}))$ is
almost-$\Fix(\si)$-internal. Let $F_i=\acl(Ea_1\ldots a_i)$ for
$i=1,\ldots,n$. By reverse induction, we may enlarge $a_n,\ldots,a_1$
so that for every $i<n$:

\begin{itemize}
\item[(a)]{$a_{i+1}$ contains a $p$-basis of $F_i(a_{i+1})_\si$ and a
transcendence basis of $F_i(a_{i+1})_\si$ over $F_i$.}

\item[(b)]{The $\si$-ideal of  difference equations satisfied by
$a_{i+1}$
over $F_i$ is generated by the difference equations satisfied by
$a_{i+1}$ over $E_i=E(a_1,\ldots,a_i)_\si$ (this is possible, since
this $\si$-ideal is finitely generated as a $\si$-ideal, see e.g.
\cite{[C]}).}
\end{itemize}
Condition (b) then implies that for every $i<n$, $E_i(a_{i+1})_\si$
and
$F_i$ are linearly disjoint over $E_i$. By the first case and (a),
$F_i(a_{i+1})_\si\subseteq F_i(a_{i+1})^s$, and the linear
disjointness of $F_i$ and $E_i(a_{i+1})_\si$ over $E_i$ then implies
that $E_i(a_{i+1})_\si\subseteq E_i(a_{i+1})^s$, so that
$E(a)_\si\subseteq E(a_1,\ldots,a_n)_\si\subseteq
E(a_1,\ldots,a_n)^s$. Then $[E(a_1,\ldots,a_n):E(a_1^p,\ldots
a_n^p)]=p^d$ where $d=tr.deg(E(a_1,\ldots,a_n)/E)<\infty$. Reasoning as
in the first case, we deduce $[E(a)_\si:E(a^p)_\si]=p^d$.

\para\vlabel{fix1}{\bf Definition of $\cals$}. Let $(K,\si)$ be an  existentially closed 
difference field of characteristic $p>0$. In this paragraph we
give a description of the classes $\cals(q)$, $q$ a non-one-based type
of SU-rank $1$.

Let $I$ be the set of pairs $(n,m)\in \nat^{>0}\times \zee$, with
$(n,m)=1$ if $m\neq 0$ and $n=1$ if $m=0$.
For each pair $(n,m)\in I$, choose a non-algebraic
type $q_{n,m}$ (over $\ffi_p^{alg}$) containing the formula
$\si^n(x^{p^m})=x$, and let
$\cals_{n,m}=\cals(q_{n,m})$. Then $q_{n,m}$ is not one-based.

By (7.1)(1) in \cite{[CHP]}, $SU(\si^n(x^{p^m})=x)=1$; as the formula
$\si^n(x^{p^m})=x$ defines a subfield of $K$, this implies that any
two
non-algebraic types containing this formula are non-orthogonal.
This observation, together with the main result of \cite{[CHP]} (see the Theorem
in section 6), shows that any type of SU-rank
$1$ which is not one-based is non-orthogonal to some $q_{n,m}$. We
define $\cals=\bigcup \cals_{n,m}$.

We will now show that if $(n,m)\neq (n',m')$ are in $I$, then
$\cals_{m,n}\cap \cals_{m',n'}=\emptyset$.

Indeed, let $F=\acl(F)$, and $a,b\in K\setminus F$ with $\si^n(a^{p^m})=a$,
$\si^{n'}(b^{p^{m'}})=b$, and assume that $a,b$ are equi-algebraic
over $F$. Then clearly $n=n'=tr.deg(F(a)_\si/F)=tr.deg(F(b)_\si/F)$.
Taking a $p^\ell$-power of $b$, we may assume that $b\in
F(a,\ldots,\si^{n-1}(a))^s$.
Let
$\tau=\si^n\Frob^m$. Then $F(a,\ldots,\si^{n-1}(a))^s$ is closed under
$\tau$ and $\tau\inv$ 
(because $\tau\si^i=\si^i\tau$ and $\tau(a)=a$), and has degree of
imperfection $n$. On the other hand, if $m\neq m'$, then the
closure under $\tau$ and $\tau\inv$ of $F(b)$ is perfect because
$\tau(b)=b^{p^{m-m'}}$. This contradicts $b\in
F(a,\ldots,\si^{n-1}(a))^s$.

\para\vlabel{thm1}{\bf Theorem}. Let $(K,\si)$ be an existentially closed  difference
field of characteristic $p>0$, let
$C\subseteq A,B$ be algebraically closed difference fields, with
$SU(B/C)<\omega$. Assume that $\acb(B/A)=A$. Then $tp(A/A\cap B)$ is
almost-$\cals$-internal.

\prf  
By Theorem \ref{prop2}, it suffices to show that whenever $A$
and $B$ satisfy the hypotheses of the theorem, then $tp(A/B)$ is \asi.
Fix such $A$, $B$, with $C=A\cap B$. We may assume $B=\acb(A/B)$; observe
that by \ref{sett}(1), 
$A=\acb(B/A)$  implies $SU(A/C)<\omega$. 

By Proposition \ref{lem37} and the discussion in \ref{fix1}, we already
know that $A=\acl(A_1\ldots A_j)$, where each $tp(A_i/C)$ is
$\cals_{n,m}$-analysable for some $(n,m)\in I$, and
$B=\acl(B_1\ldots B_j)$, where $B_i=\acb(A_i/B)$, $A_i=\acb(B_i/A)$.
If there is a counterexample to our
assertion, then
there is one where $tp(A/C)$ and $tp(B/C)$ are
$\cals_{n,m}$-analysable for some
$(n,m)\in I$, and this is what we will assume. We will also assume that
$K$ is sufficiently saturated.

Let $\tau=\si^n\Frob^m$. Let $b$ be a (finite) tuple of elements
of $B$
such that
$B=C(b)^{alg}$. Then $A$ is the smallest algebraically closed field
containing $C$ and  the field of
definition of the algebraic locus of $b$ over $A$.

We now work in the difference field $(K,\tau)$, which is a
reduct of $(K,\si)$, and is also a model of ACFA by Corollary 1.12(1)
in
\cite{[CH]}. In the reduct
$(K,\tau)$ we also have $A=\acb(Cb/A)$. By  Lemma
\ref{lem38}, we
may assume that $\tau(b)$ and $\tau\inv(b)$ are in $C(b)^s$. Hence,
there are varieties $V,W$ defined over $C$, with generics $b$ and
$(b,\tau(b))$ respectively, and with $W\subseteq V\times \tau(V)$, and
such that the
projection maps $W\to V$ and $W\to \tau(V)$ are separable
and generically finite. These maps therefore induce isomorphisms between the
jetspaces
$J^k_{(b,\tau(b))}(W)$ and $J^k_b(V)$, $J^k_{\tau(b)}(\tau(V))$ for
every $k>0$. The
proof
of Pillay and Ziegler then goes through (see chapter 3 of \cite{[PZ]}), and
shows that  $tp(A/B)$ is almost-$\Fix(\tau)$-internal (in
$(K,\tau)$). Hence there is
$M=\tau(M)^{alg}\supseteq B$, linearly disjoint  from $AB$ over $B$,
and
some tuple $a\in
\Fix(\tau)$ such that
$A\subseteq M(a)^{alg}$. Since the elements of $a$ have SU-rank $1$ in
the
difference field $(K,\tau)$, we may assume that $a$ and $A$ are
equi-algebraic over $M$.

If $n=1$, then $M=\si(M)^{alg}$, and we are done.
Assume that $n>1$; then $M$ is closed under $\si^n$ and $\si^{-n}$,
but not necessarily under $\si,\si\inv$. We need to show that there
is a difference field $(N,\si)$ extending $(B,\si)$,
containing $M$
and linearly disjoint from $AM$ over $M$, and such that
     $(N,\si^n)$ extends $(M,\tau \Frob^{-m})$. This is done
as in \cite{[CH]}, Lemma 1.12. The saturation of $(K,\si)$ then implies that
$K$
contains (a copy of) $(N,\si)$, and shows that $tp(A/B)$ is
almost-$\Fix(\tau)$-internal.

\medskip\noindent{\bf Theorem \ref{thm1}$'$}. Let $A,B$ be  difference subfields 
of $\calu$ intersecting in $C$, such that $A^{alg}\cap B^{alg}=C^{alg}$ and  $\trdeg(A/C)<\infty$. Let $D\subset B$ be generated
over $C$ by all tuples $d$ such that there exist an algebraically
closed difference field  $F$ containing $C$ and free from $B$ over $C$,
and integers $n>0$ and $m$ such that $d\in F(e)$ for some tuple $e$
of elements satisfying $\si^n\frob^m(x)=x$. Then $A$ and $B$ are free
over $D$. 

\prf When $A$ and $B$ are algebraically closed, this is a direct
consequence of \ref{thm1} and \ref{prop2}: we know that $\acb(A/B)$
realises a type over $A\cap B$ which is \asi, where $\cals$ is the
family of SU-rank $1$ types realised in some $\fix(\tau)$. Hence
$\acb(A/B)$ is contained in $D^{alg}$, which implies that $A$ and $B$
are free 
over $D$.

Assume now that $A$ and $B$ are not algebraically closed, and work over
their intersection $C$. Again, we know that $\acb(A/B)$ realises a type
over $C^{alg}$ which is \asi. Hence $\acb(A/B)$ is contained in the
maximal subset $D_0$ of $\acl(A)$ which realises an \asi-type over
$C^{alg}$. By Remark \ref{aint3}, we
have $D_0=D^{alg}$, which gives the result.

\sect{Applications of the CBP to differential and
  difference varieties}

\para {\bf Differential fields.} We will now apply some of the results of section 2 to the study of
(affine) 
differential varieties. For an introduction to the model theory of
differential fields of characteristic $0$, see e.g. \cite{Marker}. 

\medskip\noindent
{\bf Known facts}. We work in some large differentially closed field $(\calu,\delta)$ of
characteristic $0$. In analogy with the Zariski topology, we define the
{\em Kolchin topology} on each cartesian power $\calu^n$, as the topology with basic closed sets
the zero-sets of differential polynomials, which are  called {\em
 Kolchin closed sets}. This topology is Noetherian. A {\em differential
  (affine) variety} $V$ is an irreducible
Kolchin closed set.

If $A\subset \calu$ is a differential field,
then $A=\dcl(A)$ and $\acl(A)=A^{alg}$. The theory of differentially
closed fields of characteristic $0$ eliminates quantifiers and
imaginaries. 

\para
Since our results concern differential fields, we
first define the analogues of function fields and birational
morphisms. The definitions are straightforward.

\medskip
If a differential variety $V$ is defined over the differential
field $K$, we define the coordinate ring $K[V]_D$  and function field
$K(V)_D$ of $V$ as follows: let
$K[\bar X]_D$ be the ring of differential polynomials in $\bar
X=(X_1,\ldots,X_n)$, and $I$ the ideal of differential polynomials vanishing
on $V$. Then
$$K[V]_D=K[\bar X]_D/I,\qquad K(V)_D ={\rm Frac}(K[V]_D).$$
A differential variety $V$ has {\em finite order} if the transcendence
degree of $K(V)_D$ over $K$ is finite.  
If $V,W$ are differential varieties, a {\em differential-rational  map} $f:V\to W$ is
simply a map whose coordinate functions are given by elements of
$K(V)_D$; it is therefore defined on some Kolchin-open subset $U$ of
$V$. If $f(U)$ is dense in $W$ for the Kolchin topology, then we
will say that $f$ is {\em dominant}, and the map $f$ induces a $K$-embedding
of $K(W)_D$ into $K(V)_D$. Conversely, any $K$-embedding of
$K(W)_D$ into $K(V)_D$ is induced by some  dominant
differential-rational 
$f:V\to W$. A {\em finite cover}
of $V$ is a dominant differential-rational map $f:W\to V$ such that the generic fiber
of $f$ is finite. It corresponds to a finite algebraic extension $K(W)_D$ of
$K(V)_D$.  

\medskip
The {\em constant field} is 
$\calc=\{x\in\calu\mid Dx=0\}$. Any non-one-based type is 
non-orthogonal to the generic type of $\calc$. We let $\cals$ be this
generic type (over $\acl(\emptyset)$). 

\medskip
If  $A\subset \calu$, then
$K(A)_D$ denotes the differential field generated by $A$ over $K$. If $a$ is a
finite tuple, then $a$ is a {\em generic of the differential variety $V$ over
$K$} if $a\in V$ and the specialisation map $K[V]_D\to K(a)_D$ is
injective. 

\medskip
We say that a differential variety $V$ is
{\em $\calc$-internal}\footnote{Some authors say that $V$ is {\em iso-constant}.} if there is a birational $f:V\to \bar W(\calc)$  for some algebraic variety $\bar W$. We say that $V$ is
{\em almost-$\calc$-internal} if it is a finite cover of a $\calc$-internal
differential variety. This is equivalent to: if $a$ is a generic of $V$
over $K$, then $tp(a/K)$ is \asi.

\para\vlabel{aint1} {\bf Maximal almost-$\calc$-internal quotient}. Let $V$ be a
differential variety of finite 
order defined over 
the differential subfield $K$ of $ \calu$. Then $V$ has a {\em maximal
almost-$\calc$-internal quotient} $V^\#$, i.e., $V^\#$ is almost-$\calc$-internal,
and if $\pi$ is a dominant differential-rational map from $V$ to an
almost-$\calc$-internal variety $V_1$, then $\pi$ factors through $V^\#$. Furthermore, if $f:W\to V$ is a finite cover of $V$, then there is a
generically finite map $W^\#\to V^\#$.  

\prf Let $K$ be a differential field over which everything is
defined. Translated into terms of elements, this becomes: let $a$ be a
generic of $V$ over $K$, let $A=\acl(A)$ be the maximal subfield of $\acl(Ka)$ whose
type over $K$ is \asi, and let $A_0=A\cap K(a)_D$. Then
$A_0$ is finitely generated over $K$ (as a differential field or as a
field), say by a tuple $b$, and we let $V^\#$ be the differential
variety with generic $b$ over $K$, and  $V\to V^\#$ the birational map
dual to the inclusion $K(b)_D\to K(a)_D$. Note that this defines $V^\#$ uniquely
up to a differential birational correspondence, and by definition, $V^\#$
is almost-$\calc$-internal. 

Assume that $\pi:V\to V_1$ is dominant differential-rational, and let
$c=\pi(a)$. The almost-$\calc$-internality of $V_1$ is equivalent to the
almost-$\calc$-internality of $tp(c/K)$, and this implies that $c\in
A_0=K(b)_D$, and shows that the map $\pi$ factors through $V^\#$. 

For the last assertion, let $g:W\to V$ be a finite cover of $V$, let $c$
be a generic of $W$ 
such that $g(c)=a$, and let $A_1=A\cap K(c)_D$. As $c\in K(a)_D^{alg}$, we
know that $A_1$ is the maximal subfield of $K(c)_D$ which realises an
almost-$\calc$-internal type over $K$, i.e., we can take $W^\#$ to be the
differential variety of which a generator of $A_1$ over $K$ is a
generic. We clearly have $A_0\subseteq A_1\subseteq A$, and we need to show that this
extension is algebraic: but \ref{intquot} tells us that $A=A_0^{alg}$.

% . For that, it is enough to show that
% $A=A_0^{alg}$. 

% Let $d\in
% A$ be such that $A=\acl(Kd)$ and $K(d)_D\supset A_0$. Then $tp(d/K)$ is \asi, and so is the type
% of any field-conjugate of $d$ over $A_0$. It follows that the code $e$ of
% the tuple of $A_0$-conjugates of $d$ realises an \asi-type, and is
% therefore in $A_0$. As $e\in\acl(d)$, we get $A=A_0^{alg}$. 

\para\vlabel{afiber1} {\bf Maximal almost-$\calc$-internal fiber}. Let $V$ be a differential variety of finite
order defined over 
the differential subfield $K$ of $ \calu$. Then $V$ has a smallest 
quotient $V^\flat$, with generic fiber an almost-$\calc$-internal differential
variety.\footnote{In otherwords, if $\pi:V\to V_1$ is dominant with
  generic fiber almost-$\calc$-internal, then $V^\flat$ is  a quotient
  of $V_1$} Furthermore, if $f:W\to V$ is a finite cover of $V$, then there is a
generically finite map $W^\flat\to V^\flat$.  

\prf The translation in terms of differential extensions is similar to
the one done in \ref{aint1}, and
reduces the problem to the following: 

Let $B=\acl(B)$ be minimal such that $tp(a/B)$ is \asi\ (cf Theorem
\ref{thm2-1} for the existence), and let
$B_0=B\cap K(a)_D$. Then $B_0^{alg}=B$ and $tp(a/B_0)$ is \asi. But this
last statement is given by (the proof of) \ref{intfiber}.

\para\vlabel{descent1}{\bf Descent result}. For $i=1,2$, let $V_i$  be a
differential
variety of finite order defined over the differential subfield $K_i$,
 of
$\calu$, and let $k=K_1\cap K_2$. Assume that $K_1^{alg}\cap
K_2^{alg}=k^{alg}$, that $K_2$ is a regular extension of $k$, that  there is a differential rational dominant map
$f:V_1\to V_2$ defined over $(K_1K_2)^{alg}$, % that $k=K_1\cap K_2$ is
% algebraically closed,
and  that $tp(K_2/k)$ is \asi. Then there is a differential variety $V_3$ defined over
$k$, and a dominant differential rational map $g:V_2\to V_3$ such that
the generic fiber of $g$ is almost-$\calc$-internal. 

\prf Use Proposition \ref{descent} with $\dcl(\emptyset)=k$,
$b_i=K_i$, and $a_i$ a generic of $V_i$ over $K_1K_2$,
$a_2=f(a_1)$ to get $e\in K_2(a_2)_D$ such that $e\dnfo_k K_2$ and
$tp(a_2/e)$ is \asi. Since the property of \asi ity only depends on
$tp(a_2/e)$, we may take for $e$ a finite tuple. 
Our hypothesis on the extension $K_2$ of $K$
implies that $k(e)_D$ is a regular extension of $k$. If $V_3$ is the
differential locus of $e$ over $k$, and $g:V_2\to V_3$ is the dominant
map induced by the inclusion $K_2(e)_D\subset K_2(a_2)_D$, then the
generic fiber of $g$ realises an \asi-type (over $K_2(e)_D$ or $k(e)_D$). 

\para\vlabel{introdiffvar} {\bf Difference fields}. In the same vein, we now apply the results of section 2 to the study of
(affine) 
difference varieties. Again, we
 have to define the analogues of function fields and birational
morphisms. The definitions are straightforward.  

We work in some large existentially closed  difference field $\calu$. 
In analogy with the Zariski topology, we define the {\em $\si$-topology} on
each cartesian power $\calu^n$, as the topology with basic closed sets
the zero-sets of difference polynomials, which are  called {\em
$\si$-closed sets}. This topology is Noetherian. A {\em difference (affine) variety} is an irreducible
$\si$-closed set, and if this variety is defined over the difference
field $K$, we define its coordinate ring $K[V]_{\si+}$ and function field $K(V)_{\si+}$ as follows: let
$K[\bar X]_\si$ be the ring of difference polynomials in $\bar
X=(X_1,\ldots,X_n)$, and $I$ the ideal of difference polynomials vanishing
on $V$. Then
$$K[V]_{\si+}=K[\bar X]_\si/I,\qquad K(V)_{\si+} ={\rm Frac}(K[V]_{\si+}).$$ 
The {\em order} of a difference variety $V$ is the transcendence degree
of $K(V)_{\si+}$ over $K$. 
If $V,W$ are difference varieties, a {\em $\si$-rational  map} $f:V\to W$ is
simply a map whose coordinate functions are given by elements of
$K(V)_{\si+}$; it is therefore defined on some $\si$-open subset $U$ of $V$. If $f(U)$ is dense in $W$ for the $\si$-topology, then we
will say that $f$ is {\em dominant}, and the map $f$ induces a $K$-embedding
of $K(W)_{\si+}$ into $K(V)_{\si+}$. Conversely, any $K$-embedding of
$K(W)_{\si+}$ into $K(V)_{\si+}$ is induced by some %$n$ and  
dominant $\si$-rational map 
%$f:V\to W^{\si^n}$. 
$f:V\to W$. A {\em finite cover}
of $V$ is a dominant $\si$-rational map $f:W\to V$ such that the generic fiber
of $f$ is finite. It corresponds to a finite algebraic extension $K(W)_{\si+}$ of
$K(V)_{\si+}$.  

\medskip If $a$ is a tuple in $\calu$, we let $K(a)_{\si
  +}=K(\si^i(a)\mid 
i\geq 0)$; if $\si(K)=K$, then $K(a)_\si=K(\si^i(a)\mid i\in\zee)$ as in
section 3. 
We say that a tuple $a$ is a {\em generic of the difference variety
  $V$ over $K$} if
$a\in V$ and 
the natural specialisation map $K[V]_{\si +}\to K(a)_{\si +}$ is
injective. \\

We will often use  the following result (see
\cite{[C]}, 5.23.18): if $K\subset L\subset M$ are difference fields,
with $M$ finitely generated over $K$ (as a difference field), then $L$
is finitely generated over $K$.

\para {\bf Internality}. The definable closure of a difference field $K$,
$\dcl(K)$, is usually much larger than the perfect closure of $K$. The
notion of internality to $\Fix(\si)$ therefore does not have a natural
geometric intepretation. The right notion to consider is the one of
qf-internality: one replaces $\dcl$ by ``difference field generated by''. 

\smallskip\noindent
{\bf Definition}. Let $K$ be a difference field, $a$ a tuple in $\calu$,
such that $K(a)_\si/K$ is regular, let $V$ be the {\em difference locus} of $a$
over $K$ (i.e., the smallest $\si$-closed set containing $a$ and defined over
$K$), and let $\cals$ be a set of types with algebraically closed base,
which is
closed under conjugation by $\aut(\calu/K)$.  We say that $tp(a/K)$ is
{\em qf-internal to $\cals$}, or {\em qf-$\cals$-internal}, if for some
$L=\acl(L)$ containing $K$ and free from $K(a)_\si$ over
$K$, and some tuple $b$ of realisations of types in $\cals$, $a\in
L(b)_\si$. In that case, we also say that the extension $K(a)_\si/K$, and
the difference variety $V$ are {\em qf-internal to $\cals$}, or {\em
qf-$\cals$-internal}. (For the difference variety, we should really speak of
``generic'' qf-internality). And similarly we will speak of {\em \asi\
extensions}, and {\em \asi\ difference varieties}. 
Let $\tau=\si^n\Frob^m$ for some
$(m,n)\in I$ (see \ref{fix1}). If $\cals$ consists
of all types realised in  $\fix(\tau)$, then we will also 
speak of {\em qf-$\fix(\tau)$-internality}, or {\em qf-internality to $\fix(\tau)$}. 

\para{\bf Internality to fixed fields.} Let $\tau=\si^n\Frob^m$ for some
$(m,n)\in I$, and assume that $tp(a/K)$ is qf-internal
to $\Fix(\tau)$. Then one can find $L$ and $b$ as above, such that
$L(a)_\si=L(b)_\si$: take $b$ such that $L(a)_\si\cap
\fix(\tau)=(\fix(\tau)\cap L)(b)_\si$; since $L(a)_\si$ and $\fix(\tau)$
are linearly disjoint over their intersection, it follows that
$L(a)_\si$ and $L\fix(\tau)$ 
are linearly disjoint over $L(b)_\si$, and therefore
$L(a)_\si=L(b)_\si$. Note that if $m\geq 0$, then $L(b)_\si=L(b)_{\si
  +}$ and therefore also $L(a)_{\si +}=L(a)_\si$. If $m<0$, then
$L(b)_\si$ is the perfect hull of $L(b)_{\si +}$, and this implies that,
choosing $b$ so that 
 $(L\cap
\fix(\tau))(b)_{\si +}=L(a)_{\si +}\cap \fix(\tau)$, we have
$L(a)_{\si +}\supseteq L(b)_{\si +}\supseteq L(\si^j(a))_{\si +}$ for
some $j\geq 0$. If $\bar W$ is the algebraic locus of $b$ over $L$, then
there is a  purely inseparable map $\pi$ such that $\pi(V)$ is $\si$-birationally
isomorphic  (over $L$) to $\bar W(\Fix(\tau))$, the difference variety defined by $
x\in \bar W\land \tau(x)=x$.\footnote{Because $\Fix(\tau)$ is stably embedded
  (see \cite{[CHP]} 7.1), it follows that $\bar W$
is defined over $L\cap \Fix(\tau)$. If $m\geq 0$ then $j=0$ and  one does not need the
map $\pi$.} 
%%If $m\geq 0$, then $F=\fix(\tau)$ and $L(a)_{\si +}$ are linearly
%%disjoint over their intersection: assume $1=a_1,\ldots,a_n\in F$ and
%%$b_1,\ldots,b_n\in L(A)_{\si +}$ is a minimal counterexample,
%%$\sum_{i=1}^n a_ib_i=0$; then \sum_{i=1}^n \tau(a_i)\tau(b_i)=0$, and
%%as $m\geq 0$, $\tau(b_i)\in L(a)_{\si,+}$, so that we get a shorter
%%relation -- contradiction. When $m<0$, the same proof shows that $F$
%%and $L(a)_{\si+}^{1/p^\infty}$ are linearly disjoint over their
%%intersection.

\para\vlabel{facts2}{\bf Facts}. 
Let
$\tau=\si^m\frob^n$ for some $(m,n)\in I$, let $\ell\geq 1$ be an
integer. Let $K$ be a difference subfield of $\calu$, and $K'$ a
difference field isomorphic to $K$ by an isomorphism $\varphi_0$, and $\calu'$
an existentially closed difference field containing $K'$. We will work in the $\si^\ell$-difference field
$\calu[\ell]=(\calu,\si^\ell)$, and denote by $qftp(-)[\ell]$, $tp( - )[\ell]$,
$\acl_{\si^\ell}$ the quantifier-free types, types, and algebraic
closure respectively, with superscript $\calu$ or $\calu'$ if necessary.  We will use the following results:

\begin{enumerate}
\item{(\cite{[CHP]}, 1.12) Assume that $a\in K^{alg}$, and let $a'$ be a
field-conjugate 
of $a$ over $K$. Then  $qftp(a/K)[m]=qftp(a'/K)[m]$ for some
$m\geq 1$.}
\item
{(\cite{AD2}, 2.9) Let $a\in\calu$, $a'\in\calu$, and assume that there
  is an isomorphism   of $\si^\ell$-difference fields between
  $K(a)_{\si^\ell}$ and $K'(a')_{\si^\ell}$ which extends $\varphi_0$
  and sends $a$ to $a'$. Then
$tp^\calu(a/K)$ is qf-$\fix(\tau)$-internal if and only if $tp^{\calu'}(a'/K')[\ell]$ is
qf-$\fix(\tau^\ell)$-internal. }
\item(\cite{AD2}, 2.11). Let $a$ and $a'$ be as in (2).  Then
$tp^\calu(a/K)$ is one-based if and only if $tp^{\calu'}(a'/K')[\ell]$ is one-based. 
\end{enumerate}

\medskip\noindent
{\bf Conditions on the set  $\cals$}.  We fix  a set  $\cals$ of types of SU-rank $1$ with algebraically
closed base,  which is 
closed under $\aut(\calu)$-conjugation. If $p\in\cals$ is not
one-based, then for some $\tau$ as above, $p$ is non-orthogonal to any
non-algebraic type realised in $\fix(\tau)$. 
 If $\cals$ consists only of non-one-based types, then we do not
 impose any additional condition. 

If $\cals$ contains some one-based type, for convenience we will impose
% that whenever 
% $tp(a/A)\in\cals$ is one-based, and $\ell>1$, $B=\acl(B)\supset A$,
% $a\notin B$, and
% $b$ are such that $qftp(a/B)[\ell]=qftp(b/B)[\ell]$, then for any $0\leq
% i<\ell$, $tp(\si^i(b)/B)\in\cals$.
that $\cals$ contains all one-based types of SU-rank $1$. 
By abuse of language, we will
speak about 
\asi ity even when working in $\calu[\ell]$. % Our condition on $\cals$ is
% designed so that if $tp(a/K)$ is \asi, $L$ is a difference field
% extension of $K$ such that $a\in L^{alg}$, and $a'$ is a field conjugate
% of $a$ over $L$, then also $tp(a'/K)$ and $tp(\si(a')/K)$  are \asi.

\para\vlabel{aint2}{\bf Maximal \asi\ quotient of a variety}. Let $V$ be a
difference variety of finite order defined over the difference subfield $K$ of
$\calu$, and $\cals$ as above. 
Then $V$ has a maximal \asi\ quotient
$V^\#$. Furthermore, if $W$ is a finite cover of $V$, then $W^\#$ is a
finite cover of $V^\#$ via a map $\si^{-n}f$ for some integer $n$ and
tuple $f$ of rational
difference functions on $W^\#$. 

\prf Let $a$ be a generic of $V$ over $K$, and let
$A=\acl(A)\subseteq\acl(Ka)$ be maximal realising an \asi-type over
$K$. Let $A_0=A\cap K(a)_{\si +}$ and let 
$b$ be a finite tuple such that $A_0=K(b)_{\si +}$. Then $tp(b/K)$ is
\asi, which translates into: if $V^\#$ is the
difference variety of which $b$ is a generic, then $V^\#$ is
\asi, is a quotient of $V$ by a difference rational
map, and is maximal such (up to birational difference equivalence). This
is immediate observing that $K(a)_{\si +}=K(V)_{\si +}$, and
$K(b)_{\si +}=K(V^\#)_{\si +}$. 

As in the proof of \ref{aint1}, the statement about $W$ and $W^\#$
reduces to showing that $A=A_0^{alg}$. First note that because $b$
realises an \asi \ type over $K$, we have $K(b)_\si\subset A_0^{alg}$. 

As in \ref{aint1},
we argue that if $c\in A$ and $c'$ is a field conjugate of $c$ over
$K(b)_\si$, then $tp(c'/K)$ is \asi \ because $c'\in A$. Hence if $c$ is
such that $A=K(c)^{alg}$, then $c$ and $c'$ are equi-algebraic over
$K$; it then follows that the code $d$ of the  set of field conjugates of $c$
over $K(b)_\si$ is equi-algebraic with $c$ over $K$, and therefore that
$A= A_0^{alg}$: if the characteristic is $0$, then $d\in K(b)_\si$, and
if the characteristic is $p>0$, some $p^m$-power of $d$ is in $K(b)_\si$.

\para\vlabel{aint3}{\bf Remark}. A similar statement could be obtained with  maximal
qf-$\cals$-internal quotients instead: replace $A_0$ by its maximal subset $A_1$
realising a qf-$\cals$-internal type over $K$; then $A_0/A_1$ is
algebraic. \\
Observe also the following direct consequence of the proof of \ref{aint2}: let $a$ be a
tuple in $\calu$, $K$ a difference subfield of $\calu$ and $A$ the
maximal subset of $\acl(Ka)$ realising an \asi\ type over $K^{alg}$. If
$A=A_0\cap K(a)_{\si +}$, then $A=A_0^{alg}$.

\para\vlabel{afiber2}{\bf Maximal \asi\ fiber of a variety}. Let $V$ be a difference variety of finite order defined over the  difference field
$K$. Then, up to composition with a power of Frobenius, $V$ has a unique
minimal $\si$-rational quotient $V^\flat$ with the property that  
the generic fiber of the quotient map is  irreducible and
almost-$\cals$-internal. Furthermore, if $W$ is a finite cover of $V$,
then $W^\flat$ is a finite cover of $V^\flat$  via a map $f$, for some integer $n$ and
tuple $f$ of rational
difference functions on $W^\flat$.

\prf Let $a$ be a generic of $V$ over $K$, and
$A=\acl(A)\subset\acl(Ka)$ be minimal such that $tp(a/A)$ is \asi, 
let $A_0=A\cap K(a)_{\si +}$, and let  $c$ be a finite tuple such that
$A_0=K(c)_{\si +}$. We now let $V^\flat$ be the difference variety defined over
$K$ of which $c$ is a generic and $f:V\to V^\flat$ the map induced by the
inclusion $K(c)_{\si +}\subseteq K(a)_{\si +}$. 

As in the proof of \ref{afiber1}, the assertion about $W$ and $W^\flat$
reduces to showing that $A=A_0^{alg}$. 
Let $b\in A_0^s$ be such that $A=K(b)^{alg}$, and let
$b_2,\ldots,b_m$ be the 
field-conjugates of $b=b_1$ over $K(a)_\si$. By \ref{facts2}, there is $\ell\geq 1$
such that, for each $i\geq 2$, there is a $\si^\ell$-$K(a)_\si$-isomorphism
$f_i:K(a)_\si(b)_{\si^\ell}\to K(a)_\si(b_i)_{\si^\ell}$ sending $b$ to
$b_i$. Since $\si(b)\in K(b)^{alg}$, we know that
$qftp(a,\ldots,\si^{\ell-1}(a)/K(b)_{\si^\ell})[\ell]$ is
\asi, and therefore so are the types
$tp(a,\ldots,\si^{\ell-1}(a)/K(\si^j(b_i)_{\si^\ell}))[\ell]$ for $0\leq
j<\ell$ (it is clear for $j=0$; then apply powers of $\si$ to get the
result for the other values of $j$). Letting
$B=\bigcap_{i=1}^m\bigcap_{j=0}^{\ell-1} \acl_{\si^\ell}(K\si^j(b_i))$
and noting that $B=\si(B)$,  \ref{facts2} and \ref{lem2-1}
imply that
$tp(a/B)$ is \asi. The minimality of $A$ and the fact that
$b_1\in A$ now imply $A=B$. It follows that all tuples $b_i$ belong to
$A$, since 
$tr.deg(K(b_i)/K)=tr.deg(K(b)/K)$.  
Hence, if $d$ is the tuple encoding the set
$\{b_1,\ldots,b_m\}$, then $K(d)^{alg}=K(b)^{alg}$ and $tp(a/K(d)_\si)$ is
\asi. For some $n,m\geq 0$ we then have $\si^n(d^{p^m})\in K(a)_{\si +}$,
which shows  $A=A_0^{alg}$.

\para\vlabel{rem-afiber2}{\bf Remark}. The proof gives the following:
let $a$ be  a
tuple in $\calu$, $K$ a difference subfield of $\calu$ and $A$ an
algebraically closed difference subfield of $\acl(Ka)$ such that
$tp(a/A)$ is \asi. If $A_0=A\cap K(a)_{\si +}$ then
$tp(a/A_0)$ is \asi.

% \para{\bf Generating families}. Consequences can also be expressed in
% terms of generating families, following [MP]. ????

\para\vlabel{descent2} {\bf Descent of difference varieties.}
The main application of our results are given by 
Theorems  \ref{descent2} and \ref{descent3}. Theorem \ref{descent3} is
an almost optimal generalisation of 
Theorem 3.3 of \cite{AD2}. 

\medskip\noindent
{\bf Theorem}. Let $K_i$, $i=1,2$,  be difference subfields of $\calu$
with intersection $k$, and $V_i$  
difference varieties of finite order defined over $K_i$, and assume that
$k^{alg}=K_1^{alg}\cap
K_2^{alg}$. 
Assume that there is a
$\si$-rational dominant $f:V_1\to V_2$ defined over $(K_1K_2)^{alg}$, that
$tp(K_2/k)$ is \asi\ and that $K_2$ is a regular extension of $k$. Then
there is a dominant map $g:V_2\to 
V_3$, with $V_3$ a difference variety defined over $k$,  such that
the generic fiber of $g$ is  \asi.

\prf Let $a_1$ be a generic of $V_1$ over $K_1K_2$, and
$a_2=f(a_1)$. Letting $b_1=K_1$ and $b_2=K_2$, applying \ref{descent}
and using \ref{afiber2},
we obtain $e\in K_2(a_2)_\si$ such that $k(e)_\si$ and $K_2$ are
free over $k$, and $tp(a_2/k(e)_\si)$ is \asi. Moreover,  $k(e)_\si$,
being a subfield of $K_2(a_2)_\si$, is a regular extension of $k$ and is
therefore 
linearly disjoint from $K_2$ over $k$;  we may assume that $k(e)_{\si
  +}=K_2(a_2)_{\si +}\cap k(e)_\si^{alg}$. If $V_3$ is the difference
variety of which $e$ is a generic (over $K_2$), then $V_3$ is defined over $k$, and
the inclusion $K_2(e)_{\si +}\subset  K_2(a_2)_{\si +}$ gives a dominant rational difference
map $g:V_2\to V_3$
(defined over $K_2$) such that $g(a_2)=e$, and with generic fiber \asi.

\para\vlabel{descent3}{\bf Theorem}. Let $K_1,K_2$ be fields
intersecting in $k$ and with algebraic closures intersecting in $k^{alg}$; for $i=1,2$, let $V_i$ be an absolutely irreducible
variety and 
$\phi_i:V_i\to V_i$ a dominant rational map defined over $K_i$. Assume  that $K_2$ is a regular extension of $k$ and 
that there are an integer $r\geq 1$ and a dominant rational map $f:V_1\to V_2$ 
such that $f\circ \phi_1=\phi_2^{(r)} \circ f$, where $\phi_2^{(r)}$
denotes the function obtained by iterating $r$ times $\phi_2$. Then there is a variety $V_0$ and a
dominant rational map $\phi_0:V_0\to V_0$, all defined over $k$, a dominant
map $g:V_2\to V_0$ such that $g\circ \phi_2=\phi_0\circ   g$, and
$\deg(\phi_0)=\deg(\phi_2)$. 

\prf Observe that the rational map $f$ will be defined over
$(K_1K_2)^{alg}$, because this is a statement about algebraic varieties
and rational morphisms.

Let $a_1$ be a generic of $V_1$ over $K_1K_2$, and let
$a_2=f(a_1)$. Then $a_2$ is a generic of $V_2$ over $K_1K_2$. We fix an
existentially closed 
difference field $(\calu,\si)$ containing $K_2(a_2)$ and such that
$\si$ is the identity on $K_2$ and
$\si(a_2)=\phi_2(a_2)$. We fix another existentially closed field
$(\calu',\tau)$ containing $K_1K_2(a_1)$, such that $\tau$ is the
identity on $K_1K_2$, and $\tau(a_1)=\phi_1(a_1)$. Note that  $\tau$ and
$\si^r$ agree on  $K_2(a_2)$.  By abuse of notation, we let $\cals$ denote the set of
non-algebraic types of rank $1$
realised in $\fix(\si)$ when working in $\calu$, in $\fix(\tau)$ when
working in $\calu'$, and in $\Fix(\si^r)$
when working in $\calu[r]$.

Working in $\calu'$, by \ref{descent}, there is an algebraically
closed $\tau$-difference field $E$ contained in $K_2(a_2)^{alg}$ and
free from $K_2$ over 
$k$, such that $tp^{\calu',\tau}(K_2a_2/E)$ is \asi. By
\ref{rem-afiber2} and \ref{facts2}, 
we obtain  that
$tp^{\calu}(K_2a_2/E\cap K_2(a_2))[r]$ is  \asi, and therefore so is
$tp^{\calu}(K_2a_2/E)[r]$. Applying $\si^i$ for $i\geq 0$, we get that
$tp(K_2\si^i(a_2)/\si^i(E))[r]$ is \asi, and because $a_2\in
K_2(\si^i(a_2))^{alg}$ so is $tp(a_2/\si^i(E))[r]$.

Observe now that because $E=(E\cap
K_2(a_2))^{alg}$   and $\tau$ agrees with $\si^r$ on $K_2(a_2)$, we
have $\si^r(E)=E$.   
Hence, by Lemma \ref{lem2-1}, we may replace $E$ by $\bigcap_i \si^i(E)$
and assume that $\si(E)=E$. We now
reason as in \ref{afiber2} to show that if  $a_3$ is such that
$K_2(a_2)\cap E=k(a_3)$, then 
$tp(a_2/k(a_3)_\si)$ is \asi. Note that as $K_2(a_2)$ and $E$ are closed
under $\si$, so is $k(a_3)$. Hence, $\si(a_3)\in k(a_3)$. As $K_2$ is a
regular extension of $k$, and $k(a_3)\subset E$, it follows that
$k(a_3)_\si$ and $K_2$ are linearly disjoint over $k$. Letting $V_0$ be the algebraic locus of
$a_3$ over $k$, and $\phi_0$ the rational endomorphism of $V_0$ such that
$\si(a_3)=\phi_0(a_3)$, we get the desired $(V_0,\phi_0)$. The rational map $g$
is the one given by the inclusion $K_2(a_3)\subset K_2(a_2)$.

Remains the assertion about the degrees of the maps. By Lemma 1.11 of
\cite{AD2}, we have 
$1=\ld(a_2/K_2(a_3)_\si)=\ild(a_2/K_2(a_3)_\si)$, which implies 
$\deg(\phi_2)=\deg(\phi_0)$ and finishes the proof. 

\para{\bf Remarks}. \begin{itemize}
\item[(1)] As stated, the theorem says nothing when $\deg(\phi_2)=1$,
  since one can take $V_0$ of dimension $0$. 
\item [(2)]The assertion on the degrees of the map $\phi_2$ and
$\phi_0$ is weaker than the assertion that $tp(a_2/k(a_3)_\si)$ is
\asi. Note that for instance if
$c\in K_2$
is a finite tuple which generates over $k$ the
field of definition of 
$(V_2,\phi_2)$, then $c\in k(a_2)_\si$, and therefore 
$tp(c,a_2/k(a_3)_\si)$ is \asi.  This should have consequences on the
data $(V_2,\phi_2)$.  
\item[(3)] One can in fact show that the generic
  fiber of $g$ is qf-$\fix(\si)$-internal. The proof goes as follows: we
  know that there is some $a_4\in K_2(a_2)$ such that $tp(a_4/k(a_3))$
  is qf-$\fix(\si)$-internal, and $a_2\in K_2(a_4)^{alg}$. Observe that because $\ld(a_2/K_2)=\ld(a_4/K_2)$($=1$),
  the field $K_2(a_2)_\si$  is a finite extension of $K_2(a_4)_\si$. Let
  $L$ be a difference field containing $k(a_3)_\si$, linearly disjoint from
  $K_2(a_2)_\si$ over $k(a_3)_\si$, and such that $L(a_4)_\si=L(b)$ for
  some tuple $b$ in $\fix(\si)$. Enlarging $L$ if necessary, we will
  assume that $L$ is algebraically closed and that $\fix(\si)\cap L$ has
  absolute Galois group isomorphic to $\hat\zee$, so that   $\fix(\si)
  L$ contains the algebraic closure of $\fix(\si)$. It then follows by
  Lemma 4.2 of \cite{[CHS]} that $L(a_2)\subset L\fix(\si)$, which shows
  that $tp(a_2/k(a_3))$ is qf-$\fix(\si)$-internal. 
\end{itemize}
\sect{Appendix}
\para\vlabel{prop1}
{\bf Proposition}. Let $E$ and $B$ be
algebraically closed
subsets of $M$,  $b$ a tuple in $M$. Assume that $SU(B/B\cap
E)<\omega$,
that   $tp(b/B)$ is one-based, and   that $B\cap E=\acl(Bb)\cap
E$. Then $b$ is independent from $E$ over $B$.

\prf Assume  the result false, and take a counterexample with
$r=SU(B/B\cap E)-SU(B/E)$ minimal among all such $(B,E,b)$. We
may assume that 
$B\cap E=\acl(\emptyset)$, and  $E=\acb(Bb/E)$. Since
$tp(b/B)$ is one-based,
$\acl(Bb)\cap \acl(BE)\neq B$, and we may therefore assume that $b\in
\acl(BE)$.

If $r=0$, then $B\dnfo E$, so that $\acb(Bb/E)$ realises a
one-based type over $B\cap E$ (by \ref{lem1} with $\cals$ the set of one-based
types with algebraically closed base), and therefore $Bb\dnfo E$.
This
contradicts $b\in \acl(BE)$. Hence $r>0$.

Let $A=\acb(B/E)$.
We may then assume that $E$ and $b$ are equi-algebraic over $AB$: by
\ref{lem1}
     $tp(E/A)$ is one-based, and if $D=\acl(ABb)\cap E$, then $b\in
\acl(BD)$. Replace $E$ by $D$. 
Reasoning as in Step 3 of \ref{prop2}, there is $m\geq 2$, and $E$-independent realisations
$(B_1b_1),\ldots,(B_mb_m)$ of $tp(Bb/E)$ with $A\subset \acl(B_1\ldots
B_m)$ and $SU(B_m/B_1\ldots B_{m-1})>SU(B/A)$.
The induction hypothesis  implies  $b_1\dnfo_{B_1}B_2\ldots
B_m$ (since $\acl(B_1b_1)\cap \acl(B_2\ldots B_m)\subseteq
\acl(Bb)\cap E=\acl(\emptyset)$). Similarly $b_m\dnfo_{B_m}B_1\ldots
B_{m-1}$, and therefore $B_1\dnfo_{B_2\ldots B_m}b_m$.

If $E=A$, then $b\in \acl(AB)$, and
$b_1\in \acl(B_1\ldots B_m)$; by the above we get $b_1\in B_1$ which
is
absurd.

If $E\neq A$, then $E\not\subseteq \acl(B_1\ldots B_m)$ because
$E\dnfo_AB$, and   each
$b_i$ is
equi-algebraic with $E$ over $B_1\ldots
B_m$. Hence $b_1$ and $b_m$ are equialgebraic over $B_1\ldots B_m$.
However $SU(B_1/B_2\ldots B_mb_m)=SU(B_1/B_2\ldots B_m)>SU(B/E)$.
The induction hypothesis, together
with the fact that $\acl(B_1b_1)\cap \acl(B_2\ldots B_mb_m)\subseteq
\acl(B_1b_1)\cap E=\acl(\emptyset)$, gives $b_1\dnfo_{B_1}B_2\ldots
B_mb_m$, a contradiction.

\smallskip\noindent
{\bf Remark}. This result does not hold when
$SU(B/B\cap E)$
is infinite. Here is a counterexample for $T$ a completion of ACFA in
characteristic $0$. Let $a,b,c$ be generics and independent over
$\rat^{alg}$, and consider $d=ac+b$, and $e=\si(b)-b^2$. Then
$\acb(c,d/a,b)=\acl(\rat,a,b)$. Moreover, $tp(b/e)$ is one-based (by
example 6.1 of \cite{[CH]}) and has SU-rank $1$. One also has $\acl(a,b)\cap
\acl(c,d)=\rat^{alg}=\acl(\emptyset)$. Take for $(B,b,E)$ the triple 
$(\rat(a,e)^{alg}_\si,b,\rat(c,d)^{alg}_\si)$.

\para\vlabel{cor1} {\bf Proposition}. Let $tp(a/A)$ be a one-based
type
of
SU-rank $\omega^\alpha$ for some ordinal $\alpha$,
with $SU(A)<\omega$, $A=\acl(A)$,
and consider the class $\calp$ of all  types of SU-rank
$\omega^\alpha$ with
algebraically closed base, which are
non-orthogonal to $tp(a/A)$. Then $\calp$ contains a type $q$ whose base  $C$ is
contained in all bases of elements of $\calp$. If $tp(a/A)$ has
SU-rank
$1$ and is
trivial, then there is $c$ such that $SU(c/C)=1$ and $a\in \acl(Ac)$.

\prf
Assume that $tp(b/B)\in\calp$. Moving $a$, we may assume that
$a\dnfo_AB$.
By Lemma \ref{lem41}, there are
realisations
$a_1,\ldots,a_n$ of $tp(a/A)$ which are independent from $Ba$ over
$A$,  and
realisations  $b_1,\ldots,b_m$ of $tp(b/B)$ which are independent from
$A$ over $B$,  
such that $SU(a/ABa_1\ldots a_nb_1\ldots b_m)<\omega^\alpha$. Choose
such
$m,n$
minimal. Then $SU(a_1,\ldots,a_n,b_1,\ldots,b_m/AB)=\omega^\alpha(n+m)$,
and 
% the tuples $a_1,\ldots,a_n,b_1,\ldots,b_m$ are
% independent over
% $AB$, for each $i$ we have $a_i\dnfo_AB$, $b_i\dnfo_BA$, and 
$\acl(Aa_1\ldots a_n)\cap \acl(Bb_1\ldots b_m)=A\cap B=C$.
By  Proposition \ref{prop1}, we
know that
$\acl(Aaa_1\ldots a_n)\cap
\acl(Bb_1\ldots b_m)$ contains some element $d\notin C$. The usual routine
arguments then give 
$tp(d/C)\not\perp tp(a/A)$
and
$SU(d/C)=\omega^\alpha$.

Let $p_1,p_2\in\calp$, with bases $A_1,A_2$ contained in $A$. 
Because $SU(p)=\omega^\alpha$, the type $p$ has weight $1$. Hence the
inclusions $A_1,
A_2\subseteq A$  and the non-orthogonality of $p_1,p_2$ to $p$ imply
$p_1\not\perp p_2$. 
%% Then there are
%% algebraically closed sets $B_1$, $B_2$ containing $A$, tuples $a_1,a_2$
%% realising $p_1,p_2$ respectively and such that, for $i=1,2$,
%% $a_i\dnfo_{A_i}B_i$ and there is $c_i$ realising $p$ such that
%% $c_i\dnfo_AB_i$ and $a_i\dfo_{B_i}c_i$.  Conjugating by some element of
%% $\aut(M/A)$, we may assume that $B_1$ and $B_2$ are independent over
%% $A$. By the independence theorem, there is some $c$ independent from
%% $B_1B_2$ over $A$ and which realises $tp(c_1/B_1)\cup tp(c_2/B_2)$. Then for
%% $i=1,2$,
%% there are realisations $a'_i$ of $tp(a_i/B_i)$  such that
%% $a'_i\dfo_{B_1B_2}c$. Rank computations give 
%% $a'_1\dfo_{B_1B_2}a'_2$, and therefore $p_1\not\perp p_2$. 

Thus the set of bases of types in $\calp$ is
closed
under intersection, and has a smallest element, since one cannot have
an infinite decreasing sequence of algebraically closed sets of
finite SU-rank.

The last assertion follows immediately from triviality, as
non-orthogonality then implies non-almost-orthogonality.

%% \bigskip
%% \noindent
%% {\bf References}

\bigskip\noindent
UFR de Math\'ematiques,  Site Chevaleret \par\noindent
Universit\'e Paris 7 - Case 7012 \par\noindent
75205 Paris Cedex 13\par\noindent 
France \par\noindent
e-mail: {\tt zoe@math.univ-paris-diderot.fr}


\begin{thebibliography}{CHP1}
\bibitem[B]{[B]}{N. Bourbaki, XI, {\it Alg\`ebre Chap\^\i tre 5, Corps
commutatifs}, Hermann, Paris 1959.}

\bibitem[Bo]{B} E. Bouscaren, 
Proof of the Mordell-Lang conjecture for function fields, in {\em  Model
  theory and algebraic geometry},
Lecture Notes in Math., 1696, Springer, Berlin, 1998, 177 -- 196. 

% \bibitem[Bu]{Bu}{A. Buium, {\em Differential Algebraic Groups of Finite Dimension}, Lecture Notes in Math., vol
% 1506, Springer 1992.} 

\bibitem[C]{[C]}{R.M. Cohn, {\it Difference algebra}, Tracts in
Mathematics 17, Interscience Pub. 1965.}

\bibitem[CH]{[CH]}{Z. Chatzidakis, E. Hrushovski, Model theory of
difference
fields,
Trans. Amer. Math. Soc. 351 (1999),  2997 -- 3071.}

\bibitem[CH1]{AD1}{Z. Chatzidakis, E. Hrushovski, Difference fields and descent in
algebraic dynamics, I,  
Journal of the IMJ, 7 (2008) No 4,  653 -- 686.}

\bibitem[CH2]{AD2}{Z. Chatzidakis, E. Hrushovski, Difference fields and descent in
algebraic dynamics, II,  
Journal of the IMJ, 7 (2008) No 4,  687 -- 704.}


\bibitem[CHP]{[CHP]}{Z. Chatzidakis, E. Hrushovski, Y. Peterzil, Model
theory of difference fields, II: Periodic ideals and the trichotomy
in all characteristics, Proceedings of the London Math. Society (3) 85
(2002), 257  -- 311.}

\bibitem[CHS]{[CHS]}Z. Chatzidakis, C. Hardouin and M. Singer, On the Definitions of
Difference Galois Groups, in: Model Theory with
applications to algebra and analysis, I  (Z. Chatzidakis,
H.D. Macpherson, A. Pillay, A.J. Wilkie editors), Cambridge University
Press, Cambridge 2008, 73 -- 110.


\bibitem[H]{H}E. Hrushovski, private communication, November 2011.

\bibitem[HPP]{HPP}E. Hrushovski, D. Palac\'\i n, A. Pillay, On the
  Canonical Base Property, preprint July 2012, ArXiv 1205.5981. 

\bibitem[J]{J} Prerna Bihani Juhlin, Fine Structure of Dependence in
  Superstable Theories of Finite Rank, PhD thesis, University of Notre
  Dame, March 2010. 

\bibitem[M]{Marker}{D. Marker, {\em Model Theory: An Introduction}, GTM 217,
Springer-Verlag, 2002.}

\bibitem[MP]{[MP]}{R. Moosa, A. Pillay, On canonical bases and internality
criteria,  Illinois J.
Math. 52 (2008), 901 -- 917.}

\bibitem[Mo]{Mo}{R. Moosa, A Model-Theoretic Counterpart to Moishezon
    Morphisms, in {\em Models, Logics, and Higher-Dimensional Categories: A tribute to the work of Mihaly Makkai} (Eds. B. Hart, T. Kucera, A. Pillay, P. Scott, R. Seely), 
Centre de Recherches Mathematiques Proceedings \& Lecture Notes, Volume
53, American Mathematical Society 2011, 177 -- 188.}  

\bibitem[P]{P}Anand Pillay, Model-theoretic consequences of a theorem of
  Campana and Fujiki, Fundamenta 
Mathematicae, 174, (2002), 187 -- 192.


\bibitem[PW]{PW} Daniel Palacin, Frank O. Wagner, Ample thoughts, ArXiv:
  1104.0179, 2011. To appear in Notre Dame Journal of Formal Logic. 


\bibitem[PZ]{[PZ]}{A. Pillay, M. Ziegler, Jet spaces of varieties over
differential
and difference fields,  Selecta Math. (N.S.)  9  (2003),  no. 4,
579 -- 599.}

\bibitem[W1]{[W1]}{F. Wagner, {\em Simple theories}, Kluwer Academic Pub.,
Dordrecht
2000.}

\bibitem[W2]{[W2]}{F. Wagner, A note on one-basedness, J. of Symb. Logic 69
Nr 1 (2004), 34 -- 38.}
\end{thebibliography}
\end{document}